\documentclass[11pt]{article}
\usepackage{calc}
\usepackage{color}
\usepackage[english]{babel}
\usepackage{amsfonts}
\usepackage{latexsym}
\usepackage{placeins}
\usepackage{hyperref}
\let\svthefootnote\thefootnote
\newcommand\freefootnote[1]{%
	\let\thefootnote\relax%
	\footnotetext{#1}%
	\let\thefootnote\svthefootnote%
}
\usepackage{amsthm}
\usepackage{hyperref}
\usepackage{amsmath}
\usepackage{amsfonts}
\usepackage{amssymb}
\usepackage{amsthm}
\usepackage{breqn}
\usepackage{pdfpages}
\usepackage{makeidx}
\usepackage{blindtext}
\usepackage{float}
\usepackage[a4paper,top=5.5cm,bottom=4cm,left=2.5cm, right=2.5cm,foot=1cm]{geometry}
\newtheorem{teo}{Theorem}

\newtheorem{prop}[teo]{Proposition}
\newtheorem{cor}[teo]{Corollary}
\newtheorem{defi}[teo]{Definition}

\theoremstyle{remark}
\newtheorem{obss}[teo]{\bf Remark}

\usepackage{amssymb}
\usepackage{authblk}
\usepackage{amsmath}
\usepackage[cp1250]{inputenc}
\usepackage[OT4]{fontenc}

\renewcommand\Affilfont{\itshape\small}
\newcommand{\keywords}[1]{\noindent{\large{\bf Keywords:}} #1\\}
\let\LaTeXtitle\title
\renewcommand{\title}[1]{\LaTeXtitle{\Large{\textbf{#1}}}}

\title{Zero Time Discontinuity Mapping (ZDM) and Poincaré Discontinuity Mapping (PDM) in a neighbourhood of a regular grazing point of order 4 in impacting hybrid systems.}

\makeatletter
\newcommand\email[2][]%
{\newaffiltrue\let\AB@blk@and\AB@pand
	\if\relax#1\relax\def\AB@note{\AB@thenote}\else\def\AB@note{\relax}%
	\setcounter{Maxaffil}{0}\fi
	\begingroup
	\let\protect\@unexpandable@protect
	\def\thanks{\protect\thanks}\def\footnote{\protect\footnote}%
	\@temptokena=\expandafter{\AB@authors}%
	{\def\\{\protect\\\protect\Affilfont}\xdef\AB@temp{#2}}%
	\xdef\AB@authors{\the\@temptokena\AB@las\AB@au@str
		\protect\\[\affilsep]\protect\Affilfont\AB@temp}%
	\gdef\AB@las{}\gdef\AB@au@str{}%
	{\def\\{, \ignorespaces}\xdef\AB@temp{#2}}%
	\@temptokena=\expandafter{\AB@affillist}%
	\xdef\AB@affillist{\the\@temptokena \AB@affilsep
		\AB@affilnote{}\protect\Affilfont\AB@temp}%
	\endgroup
	\let\AB@affilsep\AB@affilsepx
}
\makeatother

\author[*1]{Maur\'{\i}­cio Firmino Silva Lima}
\author[2]{Tiago Rodrigo Perdig\~ao}
\affil[1]{Centro de Matemática, Computaç\~ao e Cogniç\~ao, UFABC, CEP: 09210-580, Brazil,}
\email{mauricio.lima@ufabc.edu.br}
\affil[2]{CEFET-Centro Federal tecnológico de Minas Gerais, CEP: 35.790-636, Brazil,}
\email{tiagomatt@cefetmg.br}
\begin{document}
	\pagenumbering{arabic}
	\date{}		
	\maketitle
	\freefootnote{2020 Mathematics Subject Classification:	37Cxx, 34C99, 34C25.}
	
	\freefootnote{Key words and phrases. Zero-Time Discontinuity Mapping ($ZDM$), Poincaré Discontinuity Mapping ($PDM$), Regular Grazing Point, Impacting Hybrid System }
	
	\freefootnote{The author Tiago Rodrigo Perdig\~ao is a Phd student in the postgraduate program of \textit{Centro de Matemática, Computaç\~ao e Cogniç\~ao} at UFABC.} 
	
	\freefootnote{*Corresponding author: Maur\'{\i}­cio Firmino Silva Lima.}

\maketitle
       

\begin{section}{Abstract}
In this work, we study the dynamics of a impact hybrid system. We built applications called $\textit{\textbf{ZDM}}$ (\textit{Zero Discontinuity Mapping}) and $\textit{\textbf{PDM}}$ (\textit{Poincaré Discontinuity Mapping}), for points in a neighbourhood of the 4 order regular grazing point, whose objective is to correct the behaviour of flows in the neighbourhood of this points.
\end{section}
\begin{section}{Introduction and Main Results}
	
\hspace{0.5cm}In applications, one of the most common analysed types of discontinuity-induced bifurcations (DIBs) is caused by a limit cycle of a flow becoming tangent to (i.e. grazing) with a discontinuity manifold. 

\vspace{0.3cm}
The study of such orbits is widely used by several authors, among them we highlight: in \cite{reference1}, Mats H. Fredriksson and Arne B. Nordmark investigate the study of grazing bifurcation of a stable periodic motion in a very general class of mechanical systems. In \cite{reference2}, M. Di Bernardo, C.J. Budd, A. R. Champneys, and P. Kowalczyk study of the bifurcations of grazing-type periodic orbits, in a neighbourhood of a regular grazing point of second order, for this, they use the $ZDM$ (Zero Time Discontinuity Mapping) and the $PDM$ (Poincaré Discontinuity Mapping). In addition, they illustrate the construction providing the $ZDM$ and the $PDM$ in a practical problems of 1DoF forced impact oscillators, without dissipation. In \cite{reference8} the authors analyse the different types of grazing bifurcations that may occur in a simple forced sinusoidal oscillator system in the presence of friction and a hard where the impacts occur. We also highlight recent work on dynamics with periodic orbits in \cite{reference3}, \cite{reference4}, \cite{reference5}, \cite{reference6} and \cite{reference7}

\vspace{0.3cm}
Our intention here is to generalize the studies carried out in \cite{reference2}, where the ZDM and the PDM in a regular grazing point of second order are provided, by given its expressions in a neighbourhood of a grazing regular point of order $4.$   Such maps are of great relevance to obtain information on the bifurcations of periodic orbits in impact systems, more specifically, they give conditions to find the local expression of the Poincaré map and so obtain informations on the bifurcations that may occur near a $T-$periodic hyperbolic orbit. 

In this context, we state the two main results of this paper.

\begin{teo}{(\textbf{ZDM- Zero Time Discontinuity Mapping})} \label{TeoZDM}
	Consider $(X,R)$ an impacting hybrid system, defined in $S^+ \cup \Sigma$, written in local form as in Definition \ref{localdefinition}, where $\Sigma=H^{-1} (0)$ with $0$ a regular value of $H$ and $X$ a $\mathcal{C}^3(\mathbb{R}^n)-$vector field.
	Suppose that $(X,R)$, admits at $x^*=0$ a regular grazing point of order $4.$ So, given $\mathrm{x}_1=(x_1,x_2,...,-\epsilon) \in \Pi_3=\{x \in \mathbb{R}^n;\,\, \mathcal{L}_X^3H(x)=0\},$ with $\epsilon\sim0^+,$ the $ZDM(\mathrm{x_1})$ has the form
	$$ZDM(\mathrm{x}_1)=\mathrm{x}_1-\dfrac{(4!)^{3/4}}{3!}W(\mathrm{x}_1)(\mathcal{L}_X^4H(x^*))^{1/4}\epsilon^{3/4}+ \mathcal{O}(\epsilon),$$
	where $W$ is the smooth function given by the ''reset" map $R(x)=x+W(x)\mathcal{L}_XH(x)$ from Definition \ref{localdefinition}. \newline
	\end{teo}

	\begin{teo}\label{TeoPDM}
		Consider $(X,R)$ an impacting hybrid system, defined in $S^+ \cup \Sigma$, written in local form as in Definition \ref{localdefinition}, where $\Sigma=H^{-1} (0)$ with $0$ a regular value of $H$ and $X$ a $\mathcal{C}^3(\mathbb{R}^n)-$vector field.
		Suppose that $(X,R)$, admits at $x^*=0$ a regular grazing point of order $4$, given $\mathrm{x}_1=(x_1,x_2,...,-\epsilon) \in \Pi_3=\{x \in \mathbb{R}^n;\,\, \mathcal{L}_X^3H(x)=0\}$, with $\epsilon\sim0^+$, the $PDM(\mathrm{x_1})$ has the form
		$$PDM(\mathrm{x_1})=\mathrm{x_1}-\Bigg[W(\mathrm{x_1})- \dfrac{\mathcal{L}_W\mathcal{L}_X^3H(x_1)}{\mathcal{L}_X^4H(x^*)}X(\mathrm{x_1})\Bigg]\dfrac{(4!)^{3/4}\left(\mathcal{L}_X^4H(\mathrm{x^*})\right)^{1/4}}{3!}\epsilon^{3/4}+\mathcal{O}(\epsilon),$$
			where $W$ is the smooth function given by the ''reset" map $R(x)=x+W(x)\mathcal{L}_XH(x)$ from Definition \ref{localdefinition}. \newline
		
\end{teo}

\end{section}

\section{Preliminary Results}

\hspace{0.6cm}Our main goal in this work is to generalize the maps $ZDM$ and $PDM$ to impacting hybrid system, denoted by $(X,R)$, defined in $S^+ \cup \Sigma$, having a regular grazing point of order $4$.

To this end, in which follows we will define a hybrid system and subsequently impacting hybrid system, this last one which will be our study environment throughout the work.
\begin{defi}
A hybrid system comprises a set of ODEs
$$x'=X_i(x),\,\, se\,\, x \in S_i.$$
and a set of ''reset" maps
$$R_{ij}: \Sigma_{ij}=\overline{S_i}\cap \overline{S_j}\longrightarrow \overline{S_i} \cup \overline{S_j}$$
Here, $\cup S_{i}=\mathcal{D}\subset \mathbb{R}^n$ such that each $S_i$ has non-empty interior. Each $\Sigma_{ij}=\bar{S}_i\cap\bar{S}_j$ is a $(n-1)$-dimensional manifold or an empty set. Furthermore, each $X_i$ and $R_{ij}$ are assumed to be of class $C^3(\mathbb{R}^n)$.

\begin{figure}[H]
	\centering
	\includegraphics[width=0.50\textwidth]{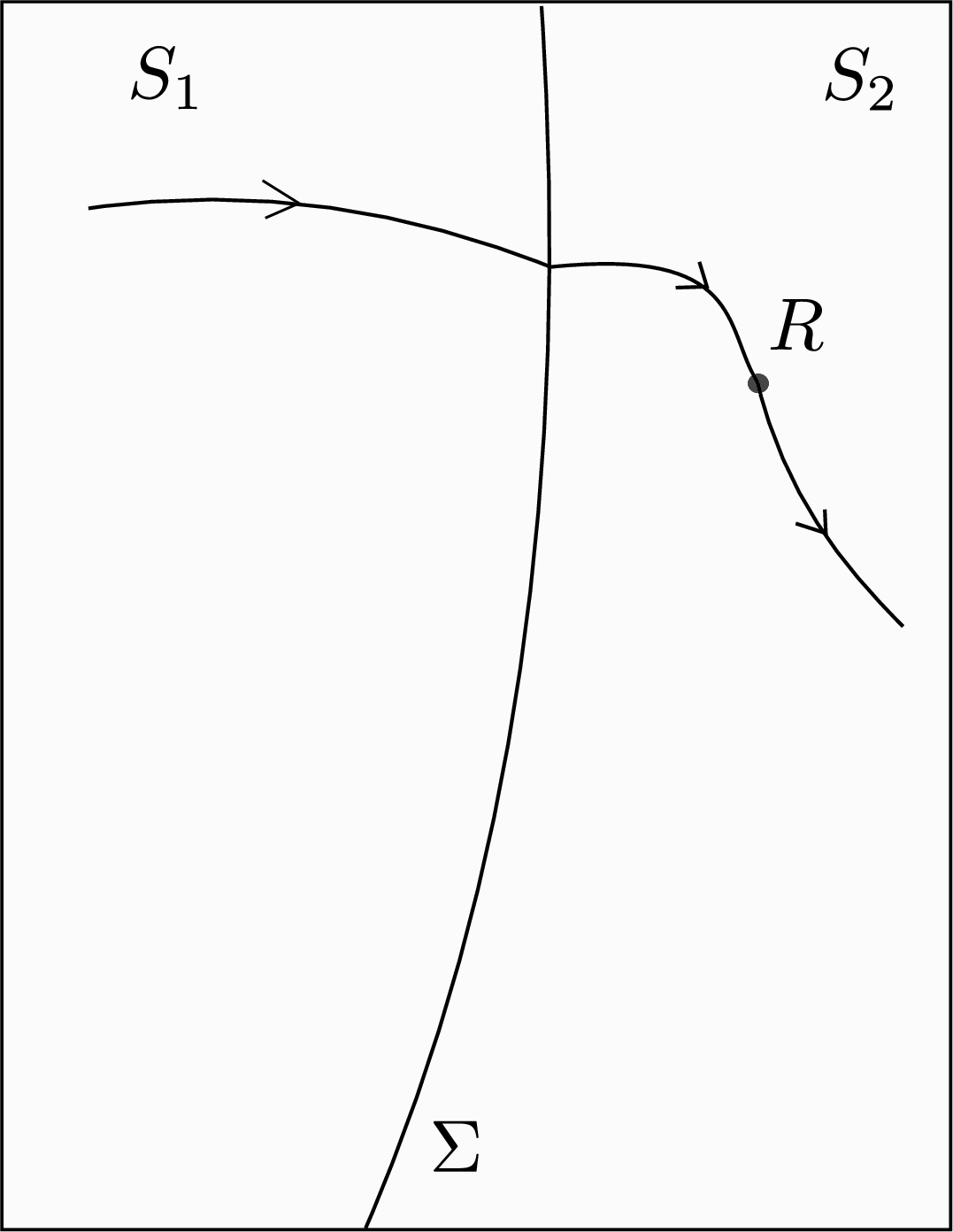}
	\caption{Hybrid System.}
\end{figure}
The subsets $\Sigma_{ij}$ are called the discontinuity boundaries or discontinuity manifolds, and each set $\Sigma_{ij}$ is a smooth submanifold of codimension 1 (i.e., locally diffeomorphic to $\mathbb{R} ^{n-1}$), locally given by $\Sigma_{ij}= \{x \in \mathbb{R}^n : H_{ij}(x) = 0\}$, for some smooth scalar function $ H_{ij} : \mathbb{R}^n \longrightarrow \mathbb{R}$, having $0$ as a regular value.
\end{defi}

Therefore we can see the discontinuity manifold $\Sigma_{ij}$ as $\Sigma_{ij}=H_{ij}^{-1}(0)$.
\begin{obss}\label{obsVar}
	Since the discontinuity manifold $\Sigma_{ij}$ is locally diffeomorphic to $\mathbb{R}^{n-1}$, then
	$$\Sigma_{ij}\cong \mathbb{R}^{n-1}\cong \{x=(x_1,...,x_n) \in \mathbb{R}^n : x_n= 0\}.$$
	In this way, if $\Sigma_{ij}= \{x \in \mathbb{R}^n : H_{ij}(x) = 0\}$, for some smooth function $H_{ij} : \mathbb {R}^n \longrightarrow \mathbb{R}$, having $0$ as a regular value then, without loss of generality, we can locally take
	$$H_{ij}(x_1,...,x_n)=x_n.$$
\end{obss}

In this work we will study a specific type of piecewise smooth system. This kind of system is object of the next definition.
\begin{defi}\label{sistemashibridodeimpacto}
An impacting hybrid system, $(X_i,R_{ij})$, is a piecewise smooth hybrid system, for which the map $R_{ij}: \Sigma_{ij}\longrightarrow \Sigma_{ij}$, and the flow is restricted to lie on one side of the boundary, that is, at $$\overline{S_i}=S_i\cup \Sigma_{ij}.$$
The map $R_{ij}$ is known as \it{''impact law"}.
\end{defi}

\hspace{0.3cm}Throughout the paper, we will consider impacting hybrid system with only one impact boundary $\Sigma$ and one impact map $R: \Sigma \longrightarrow \Sigma$, where 
\begin{equation}\label{VarDisc}
	\Sigma=H^{-1}(0)=\{x \in \mathbb{R}^n;\,\,\, H(x)=0\}, 
\end{equation}
and the vector field $X$ is defined in 
\begin{equation}\label{S^+}
S^{+}=\{x \in \mathbb{R}^n;\,\,\, H(x)>0\}.	
\end{equation}
\begin{figure}[H]
\centering
\includegraphics[width=0.50\textwidth]{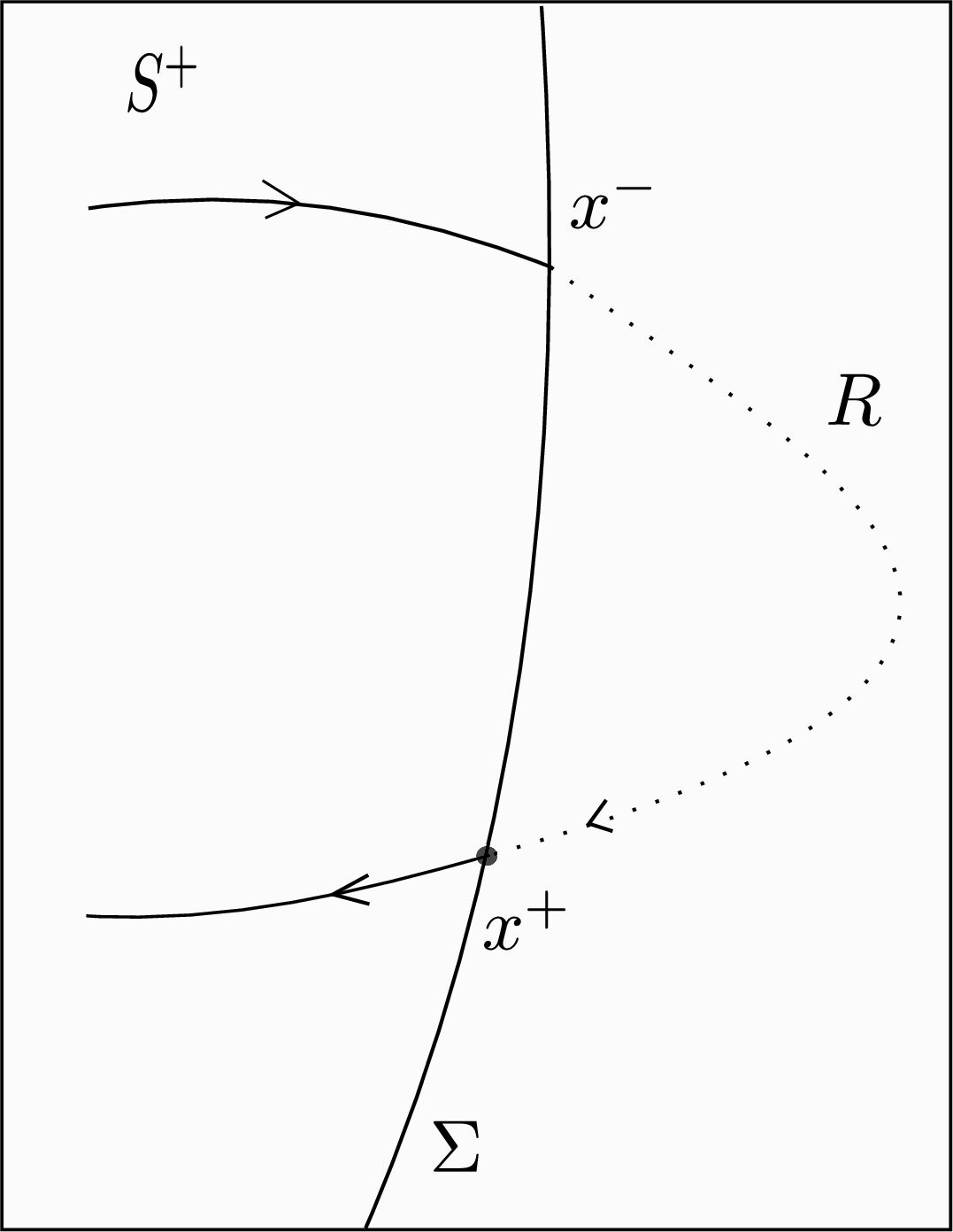}
\caption{Impacting Hybrid System.}
\end{figure}

Note that, related to the contact of the vector field $X$ with the discontinuity manifold $\Sigma$ we have two possibilities: the orbit reaches $\Sigma$ transversely at a point $x^*\in \Sigma$, or the orbit reaches $\Sigma$ tangentially at $x^*$. To distinguish these points, we define the so-called \textit{Lie Derivatives}, that measure the contact of the solutions of $X$ with the discontinuity manifold $\Sigma.$
\begin{defi}\label{Lie-derived}
	Consider $X, Y :\mathbb{R}^n \longrightarrow \mathbb{R}^n$ smooth vector fields and $H : \mathbb{R}^n \longrightarrow \mathbb{R}$ a real smooth function. The \textit{Lie Derivative $\mathcal{L}_XH(x)$} of $H$ in the direction of the vector field $X$ is defined by:
	$$\mathcal{L}_XH(x)=\dfrac{\partial H}{\partial t}(\varphi(x,t))|_{t=0}=\dfrac{\partial H}{\partial x}(x)\cdot X(x)=\nabla H(x)\cdot X(x).$$
	In a similar way, the \textit{Lie Derivative} of $\mathcal{L}_FH(x)$ in the direction of the vector field $Y$, is defined by
	$$\mathcal{L}_Y\mathcal{L}_XH(x)=\nabla\mathcal{L}_XH(x)\cdot Y(x).$$
	In the case where $Y=X$, we write $\mathcal{L}_X\mathcal{L}_XH(x)=\mathcal{L}_X^2H(x)$ and we call it the \textit{second Lie derivative} of $H $ in the direction of the vector field $X.$ More generally, we write $\mathcal{L}_X\mathcal{L}_X^{k-1}H(x)=\mathcal{L}_X^kH(x)$, for the $k-$th Lie derivative of $H$ in the direction of the vector field $X.$
\end{defi}

Lie derivatives provide the behaviour of the solutions in a neighbourhood of the discontinuity variety $\Sigma$, in particular, when the vector field is transversal or tangent at a point $x$ of $\Sigma$. As $0$ is a regular value of $H$, then for all $x\in H^{-1}(0)$, we have $\nabla H(x)\neq 0$.
\begin{figure}[H]
	\centering
	\includegraphics[width=0.80\textwidth]{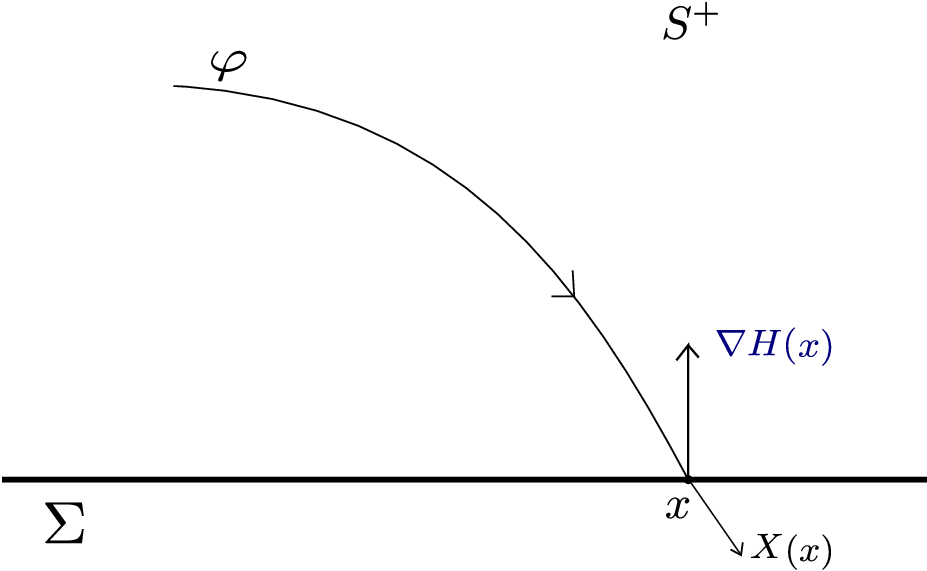}
	\caption{Geometric interpretation of the Lie derivatives of first order for piecewise smooth system.}
\end{figure}
Therefore, an orbit associated to the vector field $X$ is transversal to the discontinuity manifold at $x \in \Sigma$ if, and only if,
$$\mathcal{L}_XH(x)\neq 0.$$
On the other hand, an orbit associated to the vector field $X$ is tangent to the discontinuity manifold at $x \in \Sigma$ if, and only if,
$$\mathcal{L}_XH(x)=0.$$
Next definition characterizes impacting hybrid system studied in this paper.
\begin{defi}\label{localdefinition}
	We say that $(X,R)$, an impacting hybrid system is written in local form, if the following conditions are met
	
	1) $$
	x'=X(x), \,\,\,if\,\,\, x \in S^{+}.
	$$
	
2)	There is $W:\mathbb{R}^n\longrightarrow\mathbb{R}^n$ smooth map, such that for all $x\in \Sigma$, we have
	$$
	\begin{array}{rl}
		R(x)=&x+W(x)v,\,\,\, if\,\,\, x\in \Sigma,\\
		&\\
		v=&\mathcal{L}_XH(x)=\nabla H(x)\cdot X(x),\\
	\end{array}
	$$
	with $R: \Sigma \longrightarrow \Sigma$, $S^+=\{x \in \mathcal{D}\subset\mathbb{R}^n, H(x)>0\}$, and $W(x)$ smooth such that $x+W(x)v$ is smooth.
\end{defi}
\begin{obss}
Note that, from the above expression, $R$ acts as the identity function at tangency points of $\Sigma,$ that is,  where the vector field $X(x^*)$ is tangent to the discontinuity manifold $\Sigma$ at $x^*.$ In fact, at these points we have $\mathcal{L}_XH(x^*)=0.$ So, if $x^*\in\Sigma$ is a tangency point then,
$$R(x^*)=x^*+W(x^*)\underbrace{\mathcal{L}_XH(x^*)}_{=0}=x^*.$$
\end{obss}

\begin{figure}[H]\label{FigZDMePDM}
	\centering
	\includegraphics[width=0.9\textwidth]{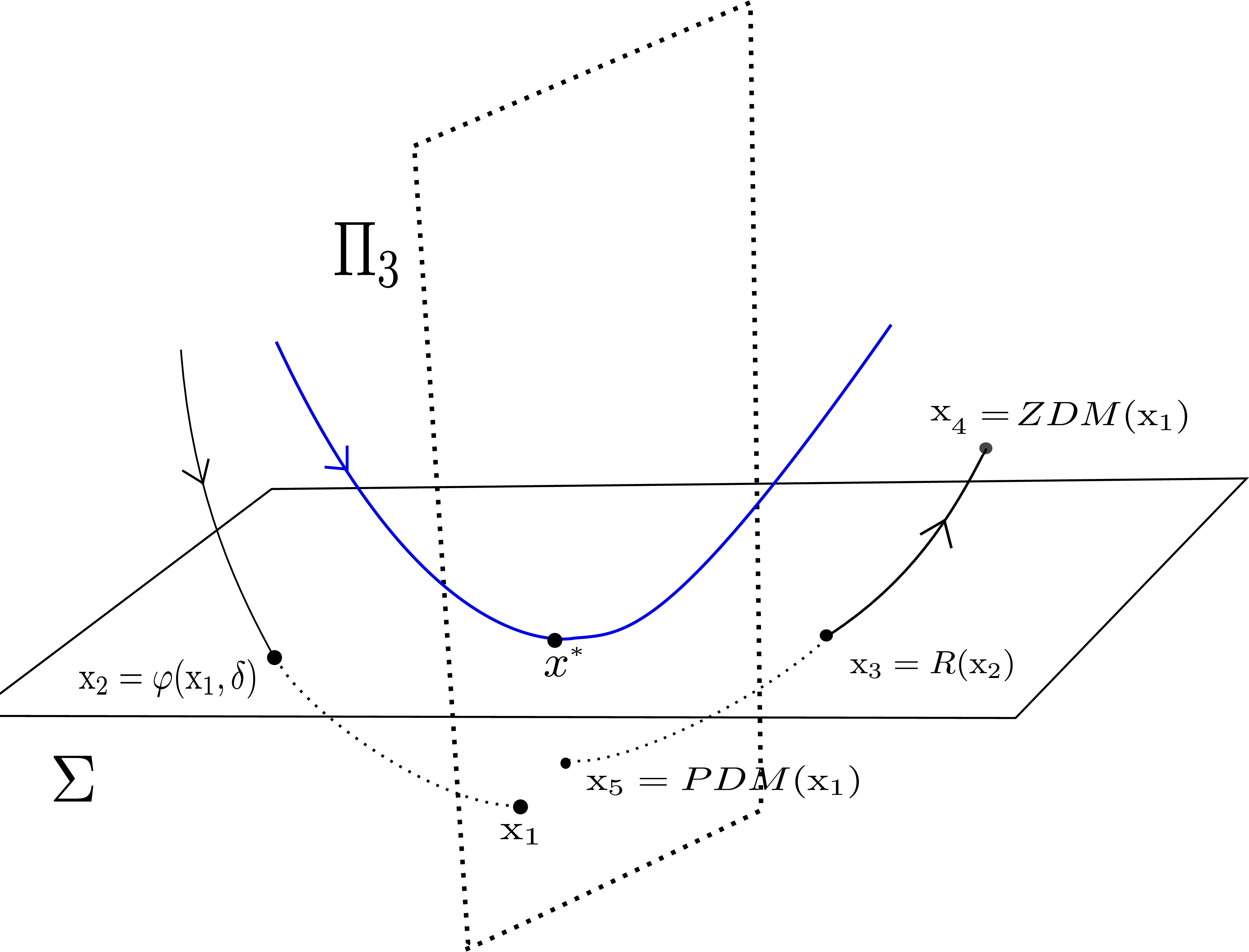}
	\caption{The geometry of the $ZDM$ and of the $PDM$ maps for a point $\mathrm{x_1} \in \Pi_3,$ satisfying $H(\mathrm{x_1})<0.$}
\end{figure}

In which follows, for points $\mathrm{x}_1 \sim x^*$, with $\mathrm{x}_1 \in \Pi_3$, the associated flow is transversal to the discontinuity manifold and passes through the point $\mathrm{x}_2=\varphi(\mathrm{x}_1,\delta)$, for $\delta<0$, (see Figure \ref{FigZDMePDM}). Therefore we have to correct the flow behaviour at the points $\mathrm{x} _1 \sim x^*$ with $H(\mathrm{x}_1)<0.$ Two maps work on this correction, namely $ZDM$ and $PDM$. The first map ($ZDM$), defined in $\Pi_3,$ is defined in such a way that the total travelling time is zero, while the second map ($PDM$), also defined in $\Pi_3,$ is given by the projection of the $ZDM$ on $\Pi_3$ under the action of the flow associated to the vector field $X.$ In order to define both maps we will take
\begin{equation}\label{Pi3}
	\Pi_3=\{x \in \mathbb{R}^n;\,\,\mathcal{L}_X^3H(x)=0\}.
\end{equation} 

Thus, for the points $\mathrm{x_1}\in\mathbb{R}^n$ satisfying $H(\mathrm{x_1})<0$ at the bottom of the discontinuity manifold, due to the discontinuity manifold $\Sigma,$ we have to correct the flow at the point $\mathrm{x_1} \in \Pi_3$. Such correction will be given by one of the following maps
$$ZDM: \mathrm{x_1}\mapsto \mathrm{x_4},$$
or
$$PDM: \mathrm{x_1}\mapsto \mathrm{x_5}.$$

The points $\mathrm{x_4}$ and $\mathrm{x_5}$ are given by the expressions
\begin{equation}\label{ZDMd}
\mathrm{x_4}=\varphi(R(\varphi(\mathrm{x_1},\delta)),-\delta),\,\, \delta<0,
\end{equation}
and
\begin{equation}
\mathrm{x_5}=\varphi(\varphi(R(\varphi(\mathrm{x_1},\delta)),-\delta),\Delta_0),\label{PDMd}
\end{equation}
in such a way that $\mathrm{x_5}$ is the projection of $\mathrm{x_4}$ on the surface $\Pi_3$. Therefore, we are able to define the $ZDM$ and the $PDM$ maps as follows:

\begin{defi}\label{defiZDM}
Consider $(X,R)$, an impacting hybrid system defined im $S^+\cup\Sigma.$ The \textit{ $ZDM$ (Zero-Time Discontinuity Mapping)} near a grazing orbit passing through the grazing point $x^*\in\Sigma,$ defined on a suitable surface $\Pi$, transversal to the flow associated to the vector field $X$ and that intersects $\Sigma$ transversely at $x^*$, is the map that does not change the evolution time of the flow which expression is given by \eqref{ZDMd}.

\end{defi}
\begin{defi}\label{defiPDM}
Consider $(X,R)$, an impacting hybrid system defined im $S^+\cup\Sigma.$ The\textit{PDM (Poincaré Discontinuity Mapping)} near a grazing orbit through the grazing point $x^*\in\Sigma,$ defined on a suitable surface $\Pi$, transverse to the flow associated to the vector field $X$ and that intersects $\Sigma$ transversely at $ x^*$, is the projection on $\Pi$ of the $ZDM(\mathrm{x_1})$ map which expression is given by \eqref{PDMd}. 
\end{defi}

In which follows, we will consider $x^* \in \Sigma$ a regular grazing point of order $4.$ This kind of point will be formalized in the next definition.

\begin{defi}\label{grazingpoint}
Given $x^* \in \Sigma$ a tangency point, we say that it is a \textit{regular grazing point of order $2k,$ $k\ge 1$}, when
$$\mathcal{L}_XH(x^*)=\mathcal{L}_X^2H(x^*)=...=\mathcal{L}_X^{2k-1}H(x^*) =0,$$
and
$$\mathcal{L}_X^{2k}H(x^*)\neq 0.$$
\end{defi}

\begin{figure}[h!]
	\centering
	\includegraphics[width=0.7\textwidth]{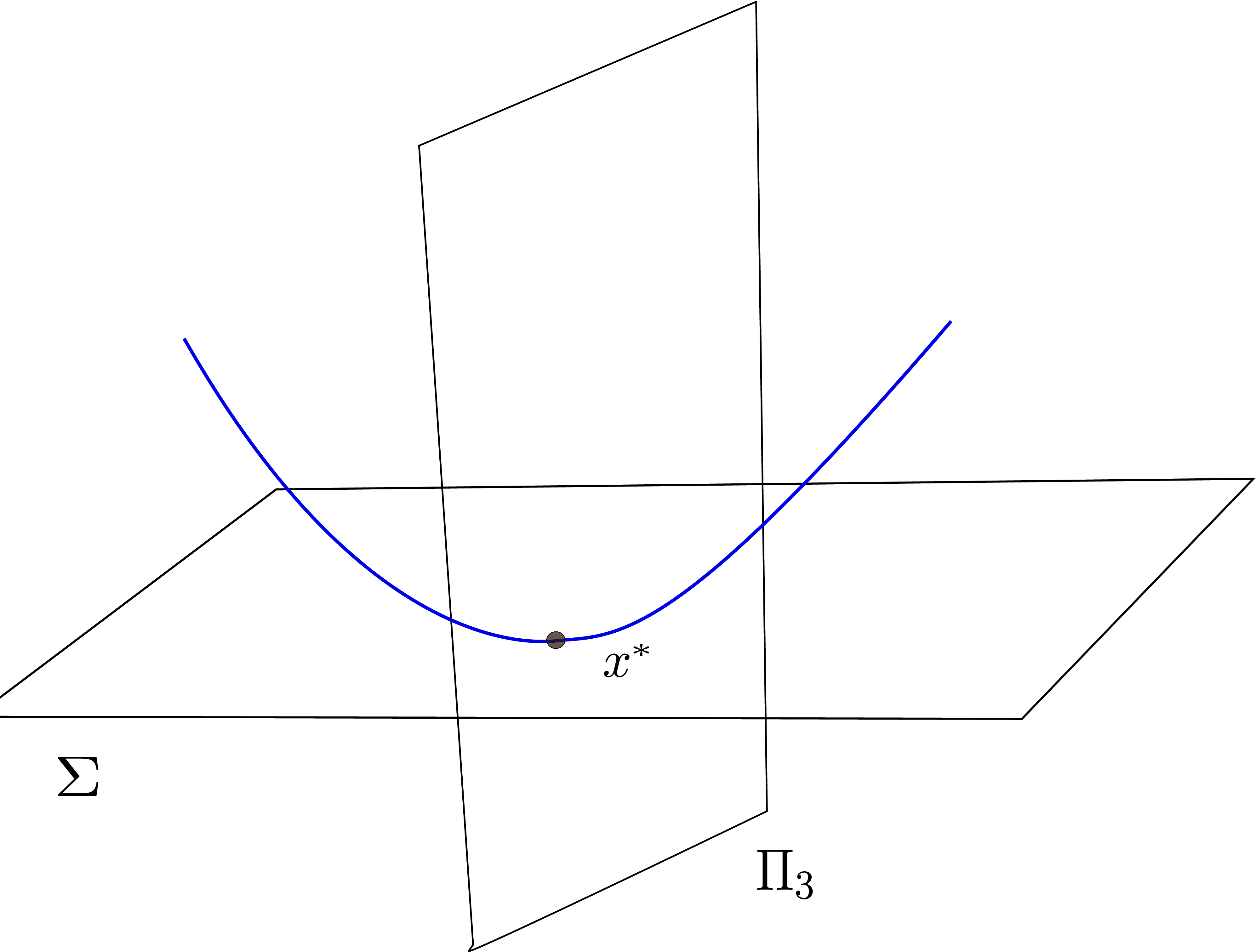}
	\caption{Grazing Regular point of order $2k,$ $k\ge 1.$}
\end{figure}

\begin{prop}Consider $x^* \in \Sigma$ a regular grazing point of order 4 of an impacting hybrid system $(X,R)$ as in Definition \ref{localdefinition}, and $\Gamma$ the orbit that passes through $x^*$. Without loss of generality, we will take as $x^*=(0,0,...,0)\in \mathbb{R}^n.$ So, related to $x^* \in \Sigma$ we have:
\begin{enumerate}
		\item $\Pi_3$ is transverse to the $(n-1)$-dimensional discontinuity manifold $\Sigma$ in $x^*$.
		\item  $\Pi_3$ is transverse to the orbit $\Gamma$ of X by $x^*$.
		
\end{enumerate}
\begin{proof}
	Note that since $x^* \in \Sigma$ is regular of order 4, then $\Pi_3=\{x \in \mathbb{R}^n; \mathcal{L}_X^3H(x)=0\}$ is such that
	$$\nabla\mathcal{L}_X^3 H(x^*)\cdot X(x^*)=\mathcal{L}_X^4 X(x^*)\neq0,$$
	that implies that $\Pi_3$ is transversal to $\Gamma$ in $x^*.$
	Moreover, as $X(x^*)$ is tangent to the discontinuity manifold $\Sigma$ at $x^*$ we also have $\Pi_3$ transversal to $\Sigma$ at $x^*\in\Sigma.$
\end{proof}	

\end{prop}

In order to get the expression for ZDM, in which follows, we will present some auxiliary results.

The first preliminary result establishes the dependence in terms of $\epsilon$ of the Lie derivatives $\mathcal{L}_X^iH, i=1,2$, at a point $\mathrm{x_1}=(x_1(\epsilon),...,-\epsilon) \in \Pi_3$, with $\Pi_3$ given in equation \eqref{Pi3}.
\begin{prop}\label{propDerivadasdeLieOrdemEpsilon}
	Consider $(X,R)$ an impacting hybrid system, defined in $S^+ \cup \Sigma$, written in local form as in Definition \ref{localdefinition}, where $\Sigma=H^{-1} (0)$, $0$ is a regular value of $H,$ and $X$ a $\mathcal{C}^3(\mathbb{R}^n)-$vector field.
	Suppose that $(X,R)$, admits at $x^*=0$ a regular grazing point of order $4.$ Given $\mathrm{x}_1=(x_1,x_2,...,-\epsilon) \in \Pi_3=\{x \in \mathbb{R}^n;\,\, \mathcal{L}_X^3H(x)=0\}$, with $\epsilon\sim0^+$, we have
	
	$$\mathcal{L}_X^iH(\mathrm{x}_1)=\mathcal{O}(\epsilon), \,\,\, i=1,2.$$
\begin{proof}
Suppose that
$$X(x_1,...,x_n)=\left(f_1(x_1,...,x_n),...,f_n(x_1,...,x_n)\right),$$ with $f_1,.., f_n \in \mathcal{C}^3(\mathbb{R}^n)$ and $H(x_1,...,x_n)=x_n$.
Given $\mathrm{x_1}=(x_1,...,-\epsilon) \in \Pi_3$, we have
$$\mathcal{L}_XH(\mathrm{x}_1)=\nabla H(\mathrm{x}_1)\cdot X(\mathrm{x}_1)=(0,0,..,1) \cdot (f_1(\mathrm{x}_1),...,f_n(\mathrm{x}_1))=f_n(\mathrm{x}_1).$$

Computing the second order Lie derivative we get
$$\mathcal{L}_X^2H(\mathrm{x}_1)=\nabla \mathcal{L}_XH(\mathrm{x_1}) \cdot X(\mathrm{x_1})=f_1(\mathrm{ x_1})\dfrac{\partial f_n}{\partial x_1}(\mathrm{x}_1)+...+f_n(\mathrm{x_1})\dfrac{\partial f_n}{\partial x_n}(\mathrm{x}_1).$$
In a similar way we obtain $\mathcal{L}_X^3H(\mathrm{x}_1)$ and $\mathcal{L}_X^4H(\mathrm{x}_1)$.

Furthermore, since $\mathrm{x_1} \sim x^*$, we can take $\mathrm{x_1}=(x_1,...,x_{n-1},-\epsilon)$ in a cube such that $x_j=\mathcal{O}(\epsilon^\alpha),$ $\alpha\ge1,$ $j=1,...,n-1.$ 

Now, expanding in Taylor series $\mathcal{L}_XH(x_1,..., -\epsilon)$ around $x^*=0$, we obtain
$$
\begin{array}{rl}
	\mathcal{L}_XH(x_1,...,x_{n-1}, -\epsilon)=&\underbrace{f_n(x^*)}_{\mathcal{L}_XH(x^*)= 0}+
	\dfrac{\partial f_n}{\partial x_1}(x^*)x_1+...+\dfrac{\partial f_n}{\partial x_n}(x^*)x_{n-1}\vspace*{0.3cm}\\
	&-\dfrac{\partial f_n}{\partial x_n}(x^*)\epsilon+\mathcal{O}(\|\mathrm{x_1}\|^2)\vspace*{0.3cm}\\
	= &\displaystyle\sum_{j=1}^n \dfrac{\partial f_n}{\partial x_j}(x^*)x_j+\mathcal{O}(\|\mathrm{x_1}\|^2)\vspace*{0.3cm}\\
	=&\displaystyle\sum_{j=1}^n \dfrac{\partial f_n}{\partial x_j}(x^*)k_j\epsilon+\mathcal{O}(\|\mathrm{x_1}\|^2).
\end{array}
$$

Similarly, we have
$$\mathcal{L}_X^{i}H(x_1,..., -\epsilon)=\mathcal{L}_X^{i}H(x^*)+\sum_{j=1}^ n\dfrac{\partial \mathcal{L}_X^{i}H}{\partial x_j}(x^*)x_j+\mathcal{O}(\epsilon^2),$$
with, $i=1,2, \,\, j=1,...,n-1$. From Definition \ref{grazingpoint} we have $\mathcal{L}_X^{i}H(x^*)=0$, $i=1,2,3$ and $\mathcal{O}(x_j)=\mathcal{O}(\epsilon^\alpha), \, \alpha \ge 1, \, j=1,...,n-1$ and $\mathcal{O}(x_n)=\mathcal{O}(\epsilon).$ Therefore the previous expression satisfy
\begin{center}
	$\mathcal{L}_X^{i}H(x_1,...,x_{n-1},-\epsilon)=\mathcal{O}(\epsilon)$,\quad$i=1,2.$
\end{center}
This concludes the result.
\end{proof}
\end{prop}

The next proposition gives the first time $\delta<0$, such that $\varphi(\mathrm{x_1},\delta) \in \Sigma$, where $\mathrm{x_1} \in \Pi_3,$ and $\varphi$ is the flow associated to the vetor field $X.$
\begin{prop}\label{propZDM}
Consider $(X,R)$ an impacting hybrid  system defined in $S^+ \cup \Sigma$, written in local form as in Definition \ref{localdefinition}, where $\Sigma=H^{-1} (0)$, $0$ is a regular value of $H$ and $X$ is a $\mathcal{C}^3(\mathbb{R}^n)-$vector field.
Suppose that $(X,R)$, admits at $x^*=0$ a regular grazing point of order $4$, given $\mathrm{x}_1=(x_1,x_2,...,-\epsilon) \in \Pi_3=\{x \in \mathbb{R}^n;\,\, \mathcal{L}_X^3H(x)=0\}$, with $\epsilon\sim0^+$, the first time $\delta<0$, such that, under the action of the flow $\varphi$, we have $\varphi(\mathrm{x}_1, \delta) \in \Sigma$, is given by
$$\delta=-\bigg(\dfrac{4!}{{\mathcal{L}_X ^4 H(\mathrm{x}_1)}}\bigg)^{1/4}\epsilon^{1 /4}+\mathcal{O}(\epsilon^{1/2}).$$
\begin{proof}
	Let $X$ be a $\mathcal{C}^3(\mathbb{R}^n)-$vector field defined in $S^+ \cup \Sigma $, where $\Sigma=H^{-1 }(0)$, with $0$ a regular value of $H.$ without loose of generality we can take
	We know that
	$$H(x_1,...,x_n)=x_n.$$
	In this case, $\Sigma=\{x=(x_1,...,x_n) \in \mathbb{R}^n;\,\, \space x_n=0\}$.
	
	Given $\mathrm{x}_1=(x_1,...,x_n) \in \Pi_{3} \cap S^-,$ $\mathrm{x_1}\sim x^*,$ we can write $\mathrm{x}_1=(x_1,...,- \epsilon)$, with $\epsilon \sim 0^+.$ So
	$$H(x_1,...,-\epsilon)=-\epsilon.$$
	
	In order to find $\delta <0$, so that $\varphi(\mathrm{x_1},\delta) \in \Sigma$, that is, $H(\varphi(\mathrm{x}_1, \delta))=0$, we will consider the following steps:
		
	Expanding the expression $H(\varphi(\mathrm{x}_1,\delta))=0$, in Taylor serie around $\delta=0$, we have
	\begin{equation}
		\begin{array}{rl}
			H(\varphi(\mathrm{x}_1,\delta))=&H(\mathrm{x}_1)+\mathcal{L}_X H(\mathrm{x}_1)\delta+\mathcal{L}_X ^2 H(\mathrm{x}_1)\dfrac{\delta^2}{2}+\mathcal{L}_X ^3 H(\mathrm{x}_1)\dfrac{\delta^3}{3 !}\\
			& +\mathcal{L}_X ^4 H(\mathrm{x}_1)\dfrac{\delta^4}{4!}+\mathcal{O}(\delta^5)\space=0.
		\end{array}
	\end{equation}
	
	So,
	\begin{equation}
		-\epsilon+\mathcal{L}_X H(\mathrm{x}_1)\delta+\mathcal{L}_X ^2 H(\mathrm{x}_1)\dfrac{\delta^2}{2}+\mathcal{L}_X ^3 H(\mathrm{x}_1)\dfrac{\delta^3}{3!}+\mathcal{L}_X ^4 H(\mathrm{x}_1)\dfrac{\delta^4}{4!}+ \mathcal{O}(\delta^5)=0.
	\end{equation}
	
	Since $\mathcal{L}_X^3H(\mathrm{x_1})=0$, we have to find $\delta<0$, such that

	\begin{equation}\label{functionK}
		K(\delta, -\epsilon)=\mathcal{L}_X H(\mathrm{x}_1)\delta+\mathcal{L}_X ^2 H(\mathrm{x}_1)\dfrac{\delta^ 2}{2}\underbrace{-\epsilon+\mathcal{L}_X^4 H(\mathrm{x}_1)\dfrac{\delta^4}{4!}+\mathcal{O}(\delta^5)}_ {g(\delta, \epsilon)}=0.
	\end{equation}
Initially, we analyse the zeroes of $g(\delta, \epsilon)=-\epsilon+\mathcal{L}_X ^4 H(\mathrm{x_1})\dfrac{\delta^4}{4!}+\mathcal{O}(\delta^5)$. By the change of variable $u=\delta^4$, we have
\begin{equation}\label{raizg}
	g(u, \epsilon)=-\epsilon+\mathcal{L}_X ^4 H(\mathrm{x}_1)\dfrac{u}{4!}+\mathcal{O}(u^{5/4})=0.
\end{equation}

Note that
$$\left\{\begin{array}{l}
	g(0,0)=0\\

	\dfrac{\partial g}{\partial u}(0,0)=\dfrac{\mathcal{L}_X ^4 H(\mathrm{x}_1)}{4!}\neq 0.
\end{array}\right.$$

Thus, it follows from the \textit{Implicit Function Theorem} that, in a sufficiently small neighbourhood of the point $(u,\epsilon)=(0,0)$, there exists a unique smooth function $u(\epsilon)$, such that $ g(u(\epsilon),\epsilon)= 0$. Furthermore,
$$u'(\epsilon)=\underbrace{u(0)}_{=0}+u'(0)\epsilon+\mathcal{O}(\epsilon^2).$$

By implicit derivation of $g(u(\epsilon),\epsilon)=0$ with respect to $\epsilon$, we obtain
$$u'(\epsilon)=-\dfrac{\dfrac{\partial g}{\partial \epsilon}(u,\epsilon)}{\dfrac{\partial g}{\partial u}(u,\epsilon)}.$$
So,
$$u'(0)=\dfrac{\dfrac{\partial g}{\partial \epsilon}(0,0)}{\dfrac{\partial g}{\partial u}(0,0)}= \dfrac{1}{\dfrac{\mathcal{L}_X ^4 H(\mathrm{x}_1)}{4!}}=\dfrac{4!}{\mathcal{L}_X ^4 H( \mathrm{x}_1)}.$$

Then,
$$u(\epsilon)=\dfrac{4!}{\mathcal{L}_X ^4 H(x_1)}\epsilon+\mathcal{O}(\epsilon^2).$$

Now, since $u=\delta^4$ it follows that
$$\overline{\delta} (\epsilon)=-\bigg(\dfrac{4!}{{\mathcal{L}_X ^4 H(\mathrm{x}_1)}}\bigg)^{1 /4}\epsilon^{1/4}+\mathcal{O}(\epsilon^{1/2}),$$ is such that $g(\bar{\delta}(\epsilon),\epsilon)=0.$

Returning to the expression in \eqref{functionK}, we see that
$$K(\delta, -\epsilon)=\mathcal{L}_X H(\mathrm{x}_1)\delta+\mathcal{L}_X ^2 H(\mathrm{x}_1)\dfrac{\delta^2}{2}+g(\delta, \epsilon)=0.$$

We assert that $K$ has a real root.

In order to find the expression of this root, we go back to \eqref{functionK}. By Proposition \ref{propDerivadasdeLieOrdemEpsilon} we know that $\mathcal {L}_X^iH(\mathrm{x_1})=O(\epsilon),\,\, i=1,2$. Take $\mathcal{L}_X^iH(\mathrm{x}_1)=\epsilon b_i+\mathcal{O}(\epsilon^2),$ with $i=1,2$. Then
$$
\begin{array}{rl}
	K(\delta, \epsilon)=&\mathcal{L}_X H(\mathrm{x}_1)\delta+\mathcal{L}_X ^2 H(\mathrm{x}_1)\dfrac{\delta^ 2}{2}+g(\delta, \epsilon)\\
	&\\
	=&\epsilon\left(\underbrace{b_1\delta+b_2\delta^2}_{f(\delta)}\right)+g(\delta,\epsilon)+\mathcal{O}(\epsilon^2).
\end{array}
$$
As $\overline{\delta} (\epsilon)=-\bigg(\dfrac{4!}{{\mathcal{L}_X ^4 H(\mathrm{x}_1)}}\bigg)^{1 /4}\epsilon^{1/4}+O(\epsilon^{1/2})$ is a root of $g$, it follows from Theorem \ref{ObsRaiz} of the Appendix
\vspace{0.3cm}
that the root of $K$ is of the form
$$
\begin{array}{rl}
	\delta=&\overline{\delta}-\epsilon\dfrac{f\left(\overline{\delta}\right)}{g_{\delta}\left(\overline{\delta}\right)}\\
	=&-\bigg(\dfrac{4!}{{\mathcal{L}_X ^4 H(\mathrm{x}_1)}}\bigg)^{1/4}\epsilon^{1/4} +O(\epsilon^{1/2})-\underbrace{\epsilon\dfrac{-b_1\bigg(\dfrac{4!}{{\mathcal{L}_X ^4 H(\mathrm{x}_1 )}}\bigg)^{1/4}\epsilon^{1/4}+O(\epsilon^{1/2})}{\frac{4}{4!}\mathcal{L}_X ^ 4 H(\mathrm{x}_1)\bigg(\dfrac{4!}{{\mathcal{L}_X ^4 H(\mathrm{x}_1)}}\bigg)^{3/4}\epsilon^{3/4}+ O(\epsilon^{5/4})}}_{O(\epsilon^{1/2})}.\\
\end{array}
$$

Therefore, $K$ has a root of the form
$$\delta(\epsilon)=-\bigg(\dfrac{4!}{{\mathcal{L}_X ^4 H(\mathrm{x}_1)}}\bigg)^{1/4}\epsilon^{1/4}+O(\epsilon^{1/2}).$$		
\end{proof}
	\end{prop}
\begin{cor}\label{CorolarioZDM}
Under the same hypotheses of Proposition \ref{propZDM}, the first time $\delta<0$, such that, under the action of the flow $\varphi$ associated to the vector field $X$, we have $\varphi(\mathrm{x }_1, \delta) \in \Sigma$, is given by
	
$$\delta=-\bigg(\dfrac{4!}{{\mathcal{L}_X ^4 H(x^*)}}\bigg)^{1/4}\epsilon^{1/4}+\mathcal{O}(\epsilon^{1/2}).$$
	\begin{proof}
By Proposition \ref{propZDM} we know that
\begin{equation}\label{delta}
	\delta=-\bigg(\dfrac{4!}{{\mathcal{L}_X ^4 H(\mathrm{x}_1)}}\bigg)^{1/4}\epsilon^{1 /4}+\mathcal{O}(\epsilon^{1/2}),
\end{equation}
for $\mathrm{x}_1 \in \Pi_3$, with $ \mathrm{x_1} \sim x^*$.
Now, as $x_j=\mathcal{O}(\epsilon^{\alpha}), \, \alpha\ge1$ and using the fact that $x^*=0$, we have
$$
\begin{array}{ll}
	\mathcal{L}_X ^4 H(\mathrm{x}_1)&=\mathcal{L}_X ^4 H(x^*)+\dfrac{\partial }{\partial x}\mathcal{L} _X ^4 H(x^*)\cdot(\mathrm{x_1} -x^*)\\
	&\\
	&=\mathcal{L}_X ^4 H(x^*)+\sum_{i=1}^n\dfrac{\partial }{\partial x_i}\mathcal{L}_X ^4 H(x^* )x_i\\
	&\\
	&=\mathcal{L}_X ^4 H(x^*)+ \mathcal{O}(\epsilon).
\end{array}
$$

So,
\begin{equation}\label{quoc}
	\begin{array}{rl}
		\dfrac{4!}{\mathcal{L}_X ^4 H(\mathrm{x}_1)}=&\dfrac{4!}{\mathcal{L}_X ^4 H(x^*)+\mathcal{O}(\epsilon)}\\
		&\\
		=&\dfrac{4!}{\mathcal{L}_X ^4 H(x^*)}\bigg(\dfrac{1}{1-(-\mathcal{O}(\epsilon))}\bigg)\\
		&\\
		=&\dfrac{4!}{\mathcal{L}_X ^4 H(x^*)}\left(1+\mathcal{O}(\epsilon)\right).
	\end{array}
\end{equation}

Returning to equation \eqref{quoc} we get
$$
\bigg(\dfrac{4!}{\mathcal{L}_X ^4 H(\mathrm{x}_1)}\bigg)^{1/4}=\left(\dfrac{4!}{\mathcal {L}_X ^4 H(x^*)}(1+\mathcal{O}(\epsilon))\right)^{1/4}
=\left(\dfrac{4!}{\mathcal{L}_X ^4H(x^*)}\right)^{1/4} \left(1+\mathcal{O}(\epsilon)\right)^{1 /4}.
$$

Expanding the term $\left(1+\mathcal{O}(\epsilon)\right)^{1/4}$ in Taylor series around $\epsilon=0$, we get $\left(1+\mathcal{O}(\epsilon)\right)^{1/4}=1+\mathcal{O}(\epsilon).$

Therefore,
$$
\begin{array}{ll}
	\left(\dfrac{4!}{\mathcal{L}_X ^4 H(\mathrm{x}_1)}\right)^{1/4}&=\left(\dfrac{4!}{\mathcal{L}_X ^4 H(x^*)}(1+\mathcal{O}(\epsilon))\right)^{1/4}\\
	&\\
	&=\left(\dfrac{4!}{\mathcal{L}_X ^4H(x^*)}\right)^{1/4} \left(1+\mathcal{O}(\epsilon)\right)^{ 1/4}\\
	&\\
	&=\left(\dfrac{4!}{\mathcal{L}_X ^4H(x^*)}\right)^{1/4}\left(1+\mathcal{O}(\epsilon)\right)\\
	&\\
	&=\left(\dfrac{4!}{\mathcal{L}_X ^4H(x^*)}\right)^{1/4}+\mathcal{O}(\epsilon).
	\end{array}
$$

Substituting the last equality in \eqref{delta} we obtain
$$\delta=-\left(\dfrac{4!}{{\mathcal{L}_X ^4 H(x^*)}}\right)^{1/4}\epsilon^{1/4} +\mathcal{O}(\epsilon^{1/2}).$$
	\end{proof}
\end{cor}

Now, we are in position to provide the expressions of the $ZDM$ and $PDM$ maps for points $\mathrm{x}_1=(x_1,x_2,...,-\epsilon) \in \Pi_3.$
\begin{section}{Proof of Theorem \ref{TeoZDM}}
	
\begin{proof}
Consider $\mathrm{x}_1 \in \Pi_3$, $H(\mathrm{x}_1)<0$ with $\Pi_3$ given by equation \eqref{Pi3}. In order to find $ZDM(\mathrm{x}_1)$ we must to consider the following steps:

\begin{enumerate}
	\item We find the first time $\delta<0$, such that, $\mathrm{x}_2=\varphi(\mathrm{x}_1,\delta) \in \Sigma.$
	\item By applying the impact map $R,$ we find $\mathrm{x}_3=R(\mathrm{x}_2).$
	\item The point $x_4 \in \mathbb{R}^n$ such that $\mathrm{x_4}=ZDM(\mathrm{x_1})$ is then obtained by
	$$\mathrm{x_4}=\varphi(R(\varphi(\mathrm{x}_1,\delta)),-\delta).$$
\end{enumerate}

It follows from Corollary \ref{CorolarioZDM}, that the first $\delta<0$ satisfying $\mathrm{x}_2=\varphi(\mathrm{x}_1,\delta)$, is given by
\begin{equation}\label{deltaNovo}
	\delta=-\bigg(\dfrac{4!}{{\mathcal{L}_X ^4 H(x^*)}}\bigg)^{1/4}\epsilon^{1/4}+\mathcal{O}(\epsilon^{1/2}).
\end{equation}

From Definition \ref{localdefinition} we have
$$R(\varphi(\mathrm{x}_1,\delta))=\varphi(\mathrm{x}_1,\delta)+W(\varphi(\mathrm{x}_1,\delta))\mathcal{L}_XH(\varphi(\mathrm{x}_1,\delta)).$$

Now we have to compute $v=\mathcal{L}_XH(\varphi(\mathrm{x}_1,\delta)).$
To do this, expanding $\mathcal{L}_XH(\varphi(\mathrm{x}_1,\delta))$ in Taylor series around $\delta=0$, we have
\begin{equation}\label{ve}
	v=\mathcal{L}_XH(\varphi(\mathrm{x}_1,\delta))=\mathcal{L}_XH(\mathrm{x}_1)+\mathcal{L}_X^2H(\mathrm {x}_1)\delta+\underbrace{\mathcal{L}_X^3H(\mathrm{x}_1)\dfrac{\delta^2}{2}}_{=0}+\mathcal{L}_X ^4H(\mathrm{x}_1)\dfrac{\delta^3}{3!}+ \mathcal{O}(\delta^4).
\end{equation}

Furthermore, it follows from Proposition \ref{propDerivadasdeLieOrdemEpsilon} that
$$\mathcal{L}_X^iH(\mathrm{x_1})=\mathcal{O}(\epsilon),\,\,i=1,2.$$

So, substituting \eqref{deltaNovo} into \eqref{ve}, we get
$$
v=-\dfrac{1}{3!}\mathcal{L}_X^4H(\mathrm{x}_1)\bigg(\dfrac{4!}{{\mathcal{L}_X^4H(x^ *)}}\bigg)^{3/4}\epsilon^{3/4}+\mathcal{O}(\epsilon).
$$

But, $\mathcal{L}_X^4H(\mathrm{x}_1)=\mathcal{L}_X^4H(x^*)+O(\epsilon)$. Therefore,
\begin{equation}\label{v}
	v=-\dfrac{(4!)^{3/4}}{3!}\left(\mathcal{L}_X^4H(x^*)\right)^{1/4}\epsilon^{ 3/4}+\mathcal{O}(\epsilon).
\end{equation}

Finally, to get
$$\mathrm{x}_4=ZDM(\mathrm{x}_1)=\varphi(\underbrace{R(\overbrace{\varphi(\mathrm{x}_1,\delta)}^{\mathrm{x }_2})}_{\mathrm{x}_3},-\delta),$$
we expand the previous expression in Taylor series around $\delta=0$,
\begin{equation}\label{x4}
	\begin{array}{ll}
		\mathrm{x}_4=&R\left(\varphi\left(\mathrm{x}_1,\delta\right)\right)-X\left(R\left(\varphi\left(\mathrm{x}_1,\delta\right)\right)\right)\delta\\
		&\\
		&+\dfrac{\partial X}{\partial x}\left(R\left(\varphi\left(\mathrm{x}_1,\delta\right)\right)\right)X\left(R\left(\varphi\left(\mathrm{x}_1,\delta\right)\right)\right)\dfrac{\delta^2}{2}\\
		&\\
		&-\Bigg(\dfrac{\partial^2 X }{\partial x^2}X^2\left(R\left(\varphi\left(\mathrm{x}_1,\delta\right)\right)\right)\\
		&\\
		&	+\bigg(\dfrac{\partial X }{\partial x}\left(R\left(\varphi\left(\mathrm{x}_1,\delta\right)\right)\right)\bigg)^2X\left(R\left(\varphi\left(\mathrm{x}_1,\delta\right)\right)\right)\Bigg)
		\dfrac{\delta^3}{3!}+\mathcal{O}(\delta^4),
	\end{array}	
\end{equation}
where, $\mathrm{x}_1 \in \Pi_{2k-1}$ and $\dfrac{\partial^2 X }{\partial x^2}X^2\left(R\left(\varphi\left(\mathrm{x}_1,\delta\right)\right)\right)$ denotes the symmetric bilinear form

$
\begin{array}{lcll}
	\dfrac{\partial^2 X }{\partial x^2}:& \mathbb{R}^{n}\times\mathbb{R}^{n}&\longrightarrow&\mathbb{R}^n\\
	& \left(X\left(R\left(\varphi\left(\mathrm{x}_1,\delta\right)\right)\right),X\left(R\left(\varphi\left(\mathrm{x}_1,\delta\right)\right)\right)\right)  &\mapsto        &	\dfrac{\partial^2 X }{\partial x^2}\left(X\left(R\left(\varphi\left(\mathrm{x}_1,\delta\right)\right)\right),X\left(R\left(\varphi\left(\mathrm{x}_1,\delta\right)\right)\right)\right),
\end{array}
$
that is,
$$\dfrac{\partial^2 X }{\partial x^2}X^2\left(R\left(\varphi\left(\mathrm{x}_1,\delta\right)\right)\right)=\dfrac{\partial^2 X }{\partial x^2}\left(X\left(R\left(\varphi\left(\mathrm{x}_1,\delta\right)\right)\right),X\left(R\left(\varphi\left(\mathrm{x}_1,\delta\right)\right)\right)\right).$$

From here we will study separately each term of the right hand of equation  \eqref{x4}. From the local form of the hybrid impact system given by Definition \ref{localdefinition} we have
\begin{equation}\label{R}
	R\left(\varphi(\mathrm{x}_1,\delta)\right)=\varphi\left(\mathrm{x}_1,\delta\right)+W\left(\varphi(\mathrm{x}_1,\delta)\right)v.	
\end{equation}

Expanding $\varphi(\mathrm{x}_1,\delta)$ in Taylor series around $\delta=0$, we obtain
\begin{equation}\label{varphi}
	\varphi(\mathrm{x}_1,\delta)=\mathrm{x}_1+X(\mathrm{x_1})\delta+\dfrac{\partial X(\mathrm{x}_1)}{\partial x}X(\mathrm{x}_1)\dfrac{\delta^2}{2!}+\Bigg(\dfrac{\partial^2X}{\partial x^2}X^2(\mathrm{x}_1)+\bigg(\dfrac{\partial X}{\partial x}\bigg)^2X^2(\mathrm{x}_1)\Bigg)\dfrac{\delta^3}{3!}+\mathcal{O}(\delta^4).
\end{equation}

So,
\begin{equation}\label{W}
	\begin{array}{rl}
		W(\varphi(\mathrm{x}_1,\delta))=&W\bigg[\mathrm{x}_1+X(\mathrm{x}_1)\delta+\dfrac{\partial X }{\partial x}(\mathrm{x}_1)X(\mathrm{x}_1)\dfrac{\delta^2}{2}\\
		&\\
		&+\Bigg(\dfrac{\partial^2 X }{\partial x^2}(\mathrm{x}_1)X^2(\mathrm{x}_1)+\bigg(\dfrac{\partial X }{\partial x}(\mathrm{x}_1)\bigg)^2X(\mathrm{x}_1)\Bigg)\dfrac{\delta^3}{3!}
		+\mathcal{O}(\delta^4)\bigg]\\
		&\\
		=&W(\mathrm{x}_1)+\dfrac{\partial W}{\partial x}(\mathrm{x}_1)X(\mathrm{x}_1)\delta+\dfrac{\partial W}{\partial x}(\mathrm{x}_1)\dfrac{\partial X }{\partial x}(\mathrm{x}_1)X(\mathrm{x}_1)\dfrac{\delta^2}{2}+\mathcal{O}(\delta^3).
	\end{array}
\end{equation}

Substituting  \eqref{varphi} and \eqref{W} into \eqref{R}, we obtain
\begin{equation}	
	\begin{array}{ll}\label{RNovo}
		R(\varphi(\mathrm{x}_1,\delta))=&\mathrm{x}_1+X(\mathrm{x}_1)\delta+\dfrac{\partial X}{\partial x}(\mathrm{x}_1)X(\mathrm{x}_1)\dfrac{\delta^2}{2}\\
		&\\
		&+\Bigg(\dfrac{\partial^2 X }{\partial x^2}(\mathrm{x}_1)X^2(\mathrm{x}_1)+\bigg(\dfrac{\partial X }{\partial x}(\mathrm{x}_1)\bigg)^2X(\mathrm{x}_1)\Bigg)\dfrac{\delta^3}{3!}\\
		&\\
		&	+\mathcal{O}(\delta^4)	+W(\mathrm{x}_1)v+\mathcal{O}(\delta v),\\
	\end{array}		
\end{equation}
where from \eqref{deltaNovo} and \eqref{v}, we have $\mathcal{O}(\delta v)=\mathcal{O}(\epsilon)$ and $\delta^4=\mathcal{O}(\epsilon)$. Now, using expansion \eqref{RNovo} and expanding $X(R(\varphi(\mathrm{x}_1,\delta)))$ in Taylor series around $\mathrm{x}_1$, we get
$$
\begin{array}{ll}
	X(R(\varphi(\mathrm{x}_1,\delta)))=&X(\mathrm{x}_1)+\dfrac{\partial X}{\partial x}(\mathrm{x}_1)X(\mathrm{x}_1)\delta+\bigg(\dfrac{\partial X}{\partial x}(\mathrm{x}_1)\bigg)^2X(\mathrm{x}_1)\dfrac{\delta^2}{2}\\
	&\\
	&+O(\delta^3)+\dfrac{\partial X}{\partial x}W(\mathrm{x}_1)v+\dfrac{\partial X}{\partial x}\dfrac{\partial W}{\partial x}(\mathrm{x}_1)X(\mathrm{x}_1)\delta v\\
	&\\
	&+\dfrac{\partial X}{\partial x}\dfrac{\partial W}{\partial x}(\mathrm{x}_1)\dfrac{\partial X }{\partial x}(\mathrm{x}_1)X(\mathrm{x}_1)\dfrac{\delta^2}{2!}v\\
	&\\
	&+\dfrac{\partial X}{\partial x}\dfrac{\partial W}{\partial x}(\mathrm{x}_1)\Bigg(\dfrac{\partial^2 X }{\partial x^2}(\mathrm{x}_1)X^2(\mathrm{x}_1)+\bigg(\dfrac{\partial X }{\partial x}(\mathrm{x}_1)\bigg)^2X(\mathrm{x}_1)\Bigg)\dfrac{\delta^3}{3!}v\\
	&\\
	&+\mathcal{O}(\delta^4,v^2).
\end{array}
$$

So,
\begin{equation}\label{-Xdelta}
		X(R(\varphi(\mathrm{x}_1,\delta))(-\delta)=	-X(\mathrm{x}_1)\delta-\dfrac{\partial X}{\partial x}(\mathrm{x}_1)X(\mathrm{x}_1)\delta^2-\bigg(\dfrac{\partial X}{\partial x}(\mathrm{x}_1)\bigg)^2X(\mathrm{x}_1)\dfrac{\delta^3}{2!}+\mathcal{O}(\delta^4,\delta v).
\end{equation}

Now, to expand $\dfrac{\partial X}{\partial x}(R(\varphi(\mathrm{x}_1, \delta)))X(R(\varphi(\mathrm{x}_1,\delta)))$ from equation \eqref{x4}, we will use expansion of $R(\varphi(\mathrm{x}_1,\delta))$ given in \eqref{RNovo}. Doing this we have
$$
\begin{array}{ll}
		\dfrac{\partial X}{\partial x}(R(\varphi(\mathrm{x}_1,\delta)))X(R(\varphi(\mathrm{x}_1,\delta)))=&\dfrac{\partial X}{\partial x}(\mathrm{x}_1)X(\mathrm{x}_1)+\bigg(\dfrac{\partial X}{\partial x}(\mathrm{x}_1)\bigg)^2X(\mathrm{x}_1)\delta\\
	&\\
	&+\dfrac{\partial^2 X}{\partial x^2}(\mathrm{x}_1)X^2(\mathrm{x}_1)\delta^2+\mathcal{O}(\delta^3),
\end{array}
$$
which implies,
\begin{equation}
	\begin{array}{ll}\label{Xdelta2}
		\dfrac{\partial X}{\partial x}(R(\varphi(\mathrm{x}_1,\delta)))X(R(\varphi(\mathrm{x}_1,\delta)))\dfrac{\delta^2}{2!}=&\dfrac{\partial X}{\partial x}(\mathrm{x}_1)X(\mathrm{x}_1)\dfrac{\delta^2}{2!}\\
		&\\
		&+\Bigg(\dfrac{\partial X}{\partial x}(\mathrm{x}_1)\Bigg)^2X(\mathrm{x}_1)\dfrac{\delta^3}{2!}+\mathcal{O}(\delta^4).
	\end{array}
\end{equation}

Now for the last term in the right hand of \eqref{x4}, we have
$$-\left(\dfrac{\partial^2 X }{\partial x^2}(R(\varphi(\mathrm{x}_1,\delta)))X^2(R(\varphi(\mathrm{x}_1,\delta)))+\left(\dfrac{\partial X }{\partial x}(R(\varphi(\mathrm{x}_1,\delta)))\right)^2X^2(R(\varphi(\mathrm{x}_1,\delta)))\right)\dfrac{\delta^3}{3!}$$
\begin{equation}\label{Xdelta3}
	=-\left(\dfrac{\partial^2 X}{\partial x^2}(\mathrm{x}_1)X^2(\mathrm{x}_1)+\left(\dfrac{\partial X }{\partial x}(\mathrm{x}_1)\right)^2X(\mathrm{x}_1)\right)\dfrac{\delta^3}{3!}+\mathcal{O}(\delta^4).
\end{equation}

Finally, substituting \eqref{Xdelta3}, \eqref{Xdelta2}, \eqref{-Xdelta} and \eqref{RNovo} into the expression of $\mathrm{x_4}$ given in \eqref{x4}, we obtain
$$
\begin{array}{rl}
	\mathrm{x}_4=&R(\varphi(\mathrm{x}_1,\delta))-X(R(\varphi(\mathrm{x}_1,\delta)))\delta\\
	&\\
	&+\dfrac{\partial X}{\partial x}(R(\varphi(\mathrm{x}_1,\delta)))X(R(\varphi(\mathrm{x}_1,\delta)))\dfrac{\delta^2}{2}\\
	&\\
	&-\Bigg(\dfrac{\partial^2 X }{\partial x^2}(R(\varphi(\mathrm{x}_1,\delta)))X^2(R(\varphi(\mathrm{x}_1,\delta)))\\
	&\\
	&+\bigg(\dfrac{\partial X }{\partial x}(R(\varphi(\mathrm{x}_1,\delta)))\bigg)^2X(R(\varphi(\mathrm{x}_1,\delta)))\Bigg)
	\dfrac{\delta^3}{3!}+\mathcal{O}(\delta^4)\\
	&\\
	=&\mathrm{x}_1+X(\mathrm{x}_1)\delta+\dfrac{\partial X }{\partial x}(\mathrm{x}_1)X(\mathrm{x}_1)\dfrac{\delta^2}{2}+
	\Bigg(\dfrac{\partial^2 X }{\partial x^2}(\mathrm{x}_1)X^2(\mathrm{x}_1)\\
	&\\
	&+\bigg(\dfrac{\partial X }{\partial x}(\mathrm{x}_1)\bigg)^2X(\mathrm{x}_1)\Bigg)\dfrac{\delta^3}{3!}
	+\mathcal{O}(\delta^4)\\
	&\\
	&+W(\mathrm{x}_1)v+O(\delta v)-X(\mathrm{x}_1)\delta\\
	&\\
	&-\dfrac{\partial X}{\partial x}(\mathrm{x}_1)X(\mathrm{x}_1)\delta^2-\bigg(\dfrac{\partial X}{\partial x}(\mathrm{x}_1)\bigg)^2X(\mathrm{x}_1)\dfrac{\delta^3}{2}+\mathcal{O}(\delta^4,\delta v)\\
	&\\
	&+\dfrac{\partial X}{\partial x}(\mathrm{x}_1)X(\mathrm{x}_1)\dfrac{\delta^2}{2!}+\Bigg(\dfrac{\partial X}{\partial x}(\mathrm{x}_1)\Bigg)^2X(\mathrm{x}_1)\dfrac{\delta^3}{2!}+\mathcal{O}(\delta^4)\\
	&\\
	&-\Bigg(\dfrac{\partial^2 X}{\partial x^2}(\mathrm{x}_1)X^2(\mathrm{x}_1)+\bigg(\dfrac{\partial X }{\partial x}(\mathrm{x}_1)\bigg)^2X(\mathrm{x}_1)\Bigg)\dfrac{\delta^3}{3!}+\mathcal{O}(\delta^4).
\end{array}
$$	

Which, after simplification, gives for the general $ZDM(x_1)$
\begin{equation}\label{ZDMequationAnterior}
	ZDM(\mathrm{ x_1})=\mathrm{x}_4=\mathrm{x}_1+W(\mathrm{x}_1)v+\mathcal{O}(\delta v).
\end{equation}

Finally, substituting the value of $v$ given by equation \eqref{v} in the previous expression and using the fact that $O(\delta v)=\mathcal{O}(\epsilon)$, we obtain
\begin{equation}\label{ZDMequation}
	ZDM(\mathrm{x}_1)=\mathrm{x_4}=\mathrm{x}_1-\dfrac{(4!)^{3/4}}{3!}W(\mathrm{x}_1)\left(\mathcal{L}_X^4H(x^*)\right)^{1/4}\epsilon^{3/4}+ \mathcal{O}(\epsilon).
\end{equation}
This concludes the result.

 \end{proof}

In the next section, the will proof Theorem \ref{TeoPDM}. this result will provide the expression of $PDM(\mathrm{x_1})$ for $\mathrm{x_1}\in\Pi_3.$
\end{section} 

\begin{section}{Proof of Theorem \ref{TeoPDM}}
	
\begin{proof}
	
\begin{figure}[H]
		\centering
		\includegraphics[width=1.1\textwidth]{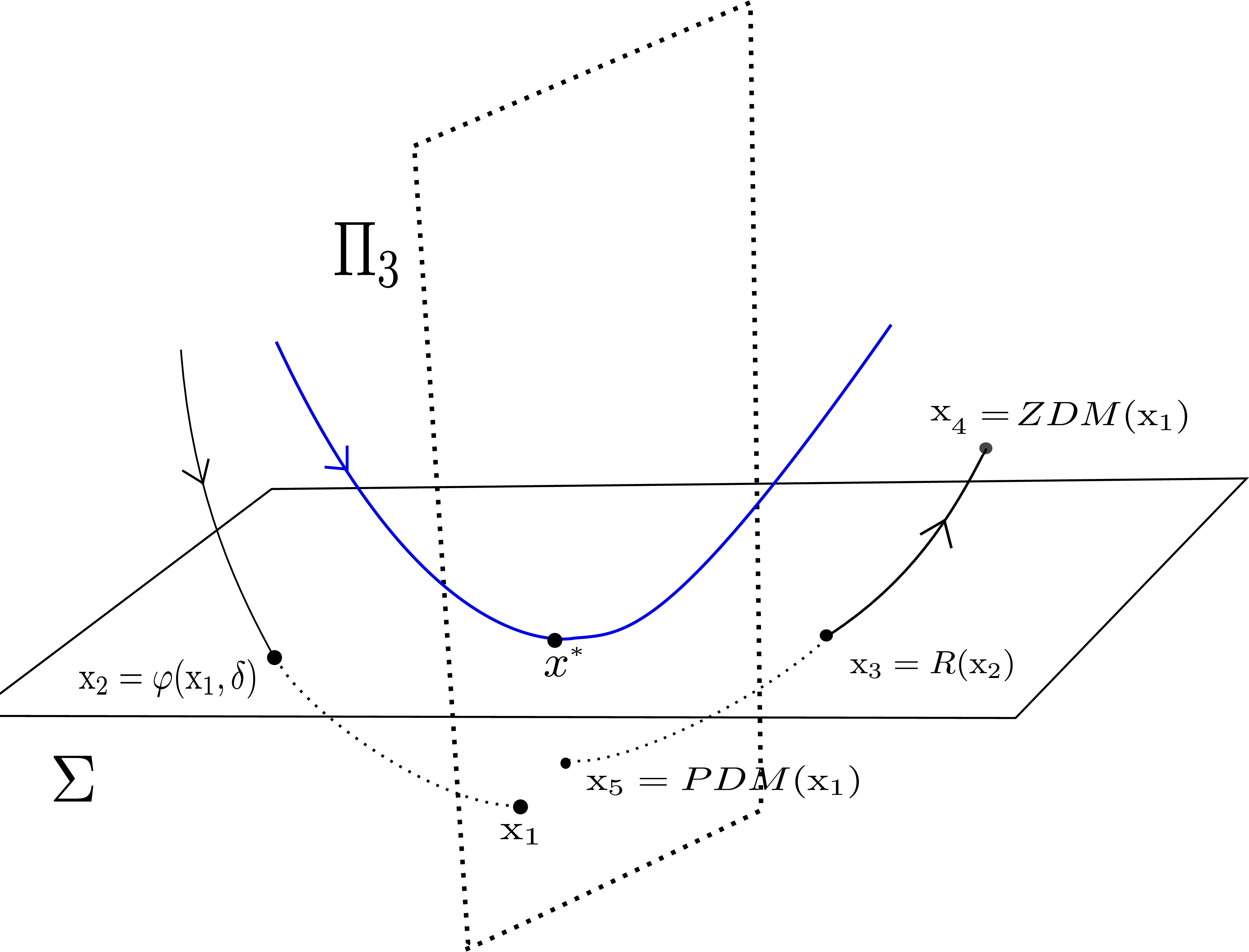}
		\caption{$PDM(\mathrm{x_1})$ for $\mathrm{x_1}\in \Pi_3.$}
	\end{figure}
The $PDM$ map is the projection of $ZDM(\mathrm{x}_1)$, obtained in Theorem \ref{TeoZDM}, on the surface $\Pi_3=\{x \in D:\,\, \mathcal{L}_X^3H(x)=0 \}.$

In this way, the $PDM(\mathrm{x_1})$ is given by the point $\mathrm{x_5}$, where
$$\mathrm{x_5}=\varphi(\mathrm{x_4},\Delta_0),$$ with $\mathcal{L}_X^3H(\mathrm{x_5})=0.$

Thus, in order to find the expression for $\mathrm{x_5}$, we will follow the steps:
\begin{enumerate}
	\item Find the first time $\Delta_0$ such that $\mathrm{x_5}=\varphi(\mathrm{x_4},\Delta_0) \in \Pi_3$, that is, $\mathcal{L}_X^3H(\mathrm{x_5})=0$.
	\item Expand the expression $\mathrm{x_5}=\varphi(\mathrm{x_4},\Delta_0)$ in Taylor series around $\Delta_0=0$.
\end{enumerate}

Therefore, initially we have to find $\Delta_0$ such that $\mathcal{L}_X^3H(\mathrm{x_5})=0$. Expanding the expression $\mathcal{L}_X^3H(\mathrm{x_5})=0$ in Taylor serie around $\Delta_0=0$, we get
\begin{equation}\label{L3}
	\mathcal{L}_X^3H(\mathrm{x_5})=\mathcal{L}_X^3H\left(\varphi\left(\mathrm{x_4},\Delta_0\right)\right)=\mathcal{ L}_X^3H(\mathrm{x_4})+\mathcal{L}_X^4H(\mathrm{x_4})\Delta_0+\mathcal{L}_X^5H(\mathrm{x_4})\dfrac{\Delta_0 ^2}{2}+\mathcal{O}\left(\Delta_0^3\right)=0.
\end{equation}
Expanding the term $\mathcal{L}_X^3H(\mathrm{x_4})$ around $\mathrm{x}_1$, with $\mathrm{x}_4$ given by equation \eqref{ZDMequationAnterior} we have
\begin{equation}
	\begin{array}{ll}\label{L3New}
		\mathcal{L}_X^{3}H(\mathrm{x_4})&=\mathcal{L}_X^3H\left(\mathrm{x_1}+W(\mathrm{x_1})v+\mathcal{O}(\delta v)\right)\\
		&\\
		&=\mathcal{L}_X^{3}H(\mathrm{x}_1)+\mathcal{L}_W\mathcal{L}_X^{3}H(\mathrm{x}_1)v+\mathcal{O}(\delta v).\\
	\end{array}
\end{equation}

Proceeding in the same way with the term $\mathcal{L}_X^4H(\mathrm{x_4}),$ we obtain
\begin{equation}\label{L4}
	\mathcal{L}_X^4H(\mathrm{x_4})=\mathcal{L}_X^4H(\mathrm{x_1})+\mathcal{O}(v).
\end{equation}

Now substituting \eqref{L3New} and \eqref{L4} in \eqref{L3} we have
$$
	\begin{array}{ll}
		\mathcal{L}_X^{3}H(\mathrm{x_5})&=\underbrace{\mathcal{L}_X^{3}H(\mathrm{x_1})}_{=0}+\mathcal{L}_W\mathcal{L}_X^{3}H(\mathrm{x_1})v+\mathcal{L}_X^4H(\mathrm{x_1})\Delta_0+\mathcal{O}(\Delta_0v, v^2,\Delta_0^2)\\
		&\\
		&=\mathcal{L}_W\mathcal{L}_X^{3}H(\mathrm{x_1})v+\mathcal{L}_X^{4}H(\mathrm{x_1})\Delta_0+\mathcal{O}(\delta v,\Delta_0 v, v^2,\Delta_0^2)=0
	\end{array}
$$

Consider
 $$G(\mathrm{x_1},v,\Delta_0)=\mathcal{L}_W\mathcal{L}_X^3H(\mathrm{x_1})v+\mathcal{L}_X^4H(\mathrm{x_1})\Delta_0+\mathcal{O}(\delta v,\Delta_0 v, v^2,\Delta_0^2)=0.$$
Note that
	$$\left\{\begin{array}{ll}
		G(x^*,0,0)=0,&\\
		&\\
		\dfrac{\partial G}{\partial \Delta_0}(0,0)= \mathcal{L}_X^4H(x^*)\neq 0.\\
	\end{array}\right.$$ 
	
It follows from the \textit{Implicit Function Theorem} that, for $(\mathrm{x_1},v) \sim (x^*,0)$ there is a unique smooth function $\Delta_0(\mathrm{x_1},v)$
such that $G(\mathrm{x_1},v,\Delta_0(\mathrm{x_1},v))=0$ and $\Delta_0(x^*,0)=0$.

	A direct computation shows that
	$$\Delta_0=-\dfrac{\mathcal{L}_W\mathcal{L}_X^3H(\mathrm{x_1})}{\mathcal{L}_X^4H(x^*)}v+\mathcal{O}(v^2).$$
	
Moreover, from equation \eqref{v} we have
	$$v=-\dfrac{(4!)^{3/4}}{3!}(\mathcal{L}_X^4H(x^*))^{1/4}\epsilon^{3/4}+\mathcal{O}(\epsilon).$$
	
So,
\begin{equation}\label{Delta_0}
	\Delta_0=\dfrac{(4!)^{3/4}}{3!}\dfrac{\mathcal{L}_W\mathcal{L}_X^3H(\mathrm{x_1})}{\mathcal{L}_X^4H(x^*)}(\mathcal{L}_X^4H(x^*))^{1/4}\epsilon^{3/4}+\mathcal{O}(\epsilon).
\end{equation}

Now we expand the expression
	 $\mathrm{x_5}=\varphi(\mathrm{x_4},\Delta_0)$ in Taylor series around $\Delta_0=0$ to obtain
	 \begin{equation}\label{PDManterior}
 	\mathrm{x_5}=\varphi(\mathrm{x_4},0)+X(\mathrm{x_4})\Delta_0+\mathcal{O}(\Delta_0^2).
	 \end{equation}

Expanding $X(\mathrm{x_4})$ around $\mathrm{x_1}$, with $\mathrm{x_4}$ given in \eqref{ZDMequationAnterior}, we have
	$$X(\mathrm{x_4})=X(\mathrm{x_1})+\dfrac{\partial X}{\partial x}W(\mathrm{x_1})v+\mathcal{O}(\delta v).$$
So, the expression of $PDM(\mathrm{x_4})$ becomes
\begin{equation}\label{PDM1}
	PDM(\mathrm{x_1})=\mathrm{x_5}=\underbrace{\mathrm{x_4}}_{ZDM(\mathrm{x_1})}+X(\mathrm{x_1})\Delta_0+ \dfrac{\partial X}{\partial x}W(\mathrm{x_1})v\Delta_0 + \mathcal{O}(\delta v, \Delta_0^2).
\end{equation}
	
Finally, replacing

	$$ZDM(\mathrm{x_1})=\mathrm{x_1}-\dfrac{(4!)^{3/4}}{3!}W(\mathrm{x_1})(\mathcal{L}_X^4H(x^*))^{1/4}\epsilon^{3/4}+ \mathcal{O}(\epsilon)$$
	and 
	$$\Delta_0=\dfrac{(4!)^{3/4}}{3!}\dfrac{\mathcal{L}_W\mathcal{L}_X^3H(\mathrm{x_1})}{\mathcal{L}_X^4H(x^*)}(\mathcal{L}_X^4H(x^*))^{1/4}\epsilon^{3/4}+\mathcal{O}(\epsilon),$$
	in the expression para a $PDM(\mathrm{x_1})$ given in \eqref{PDM1}, gives
	$$
\begin{array}{ll}
PDM(\mathrm{x_1})=&\mathrm{x_1}-\dfrac{(4!)^{3/4}}{3!}W(\mathrm{x_1})\left(\mathcal{L}_X^4H(x^*)\right)^{1/4}\epsilon^{3/4}\\
&\\
		&+\dfrac{(4!)^{3/4}}{3!}X(\mathrm{x_1})\dfrac{\mathcal{L}_W\mathcal{L}_X^3H(\mathrm{x_1})}{\mathcal{L}_X^4H(x^*)}\left(\mathcal{L}_X^4H(x^*)\right)^{1/4}\epsilon^{3/4}+\mathcal{O}(\epsilon),
	\end{array}
	$$

or
\begin{equation}\label{PDMequation}
	PDM(\mathrm{x_1})=\mathrm{x_1}-\Bigg[W(\mathrm{x_1})- \dfrac{\mathcal{L}_W\mathcal{L}_X^3H(x_1)}{\mathcal{L}_X^4H(x^*)}X(\mathrm{x_1})\Bigg]\dfrac{(4!)^{3/4}\left(\mathcal{L}_X^4H(x^*)\right)^{1/4}}{3!}\epsilon^{3/4}+\mathcal{O}(\epsilon).
\end{equation}
This concludes the proof.

\end{proof}

\section{Example: ZDM and PDM for an Impacting Hybrid System defined by a Perturbed Hamiltonian System}

\medskip

Consider the field $X(x,y)$ given by the differential system

\begin{equation}\label{HamilPert}
	\begin{array}{ll}
		x^\prime &= -(y-1)^3=f(x,y),\\
		&\\
		y^\prime &=x^3-\xi(x^4+(y-1)^4-1)=g(x,y),
	\end{array}
\end{equation}
and the impacting hybrid system $(X,R)$ defined in $\cup S^+\cup\Sigma$ with $\Sigma=\{(x,y);\,H(x,y)=y=0\},$ $S^+=\{(x,y);\,H(x,y)>0\}$ and, by Definition \ref{localdefinition}, $R(x,0)=(x,0)+W(x,0)\mathcal{L}_XH(x,0)$ where $W(x,y)$ is a smooth function defined in $\mathbb{R}^2$ such that for all $(x,0)\in \Sigma$ we have $W(x,0)=(w_1(x),0).$\\

Note that for $\xi=0$ differential system \eqref{HamilPert} becomes
\begin{equation}\label{Hamil}
	\begin{array}{ll}
		x^\prime &= -(y-1)^3,\\
		y^\prime &=x^3.
	\end{array}
\end{equation}
Differential system \eqref{Hamil} is a Hamiltonian system, with Hamiltonian function
$$\tilde{H}(x,y)=-\dfrac{x^4}{4}-\dfrac{(y-1)^4}{4}.$$
Moreover, the perturbation $\xi$ preserves $\Gamma: x^4+(y-1)^4=1$ as a periodic orbit of \eqref{HamilPert}.

Consider $x^*=(0,0)\in\Gamma.$ We will see that $x^*$ is a regular grazing point of order $4$ of the impacting hybrid system $(X,R),$ then we will compute the $ZDM$ and $PDM$ maps, for points $\mathrm{x_1} \in \Pi_3$ and $\mathrm{x_1}\sim x^*.$
Observe that
$$
\begin{array}{ll}
	\mathcal{L}_XH(x,y)&=\nabla H(x,y)\cdot X(x,y)\\
	&\\
	&=(0,1)\cdot(-(y-1)^3,x^3-\xi(x^4+(y-1)^4-1))\\
	&\\
	&=x^3-\xi(x^4+(y-1)^4-1=g(x,y).
\end{array}
$$

So,
$$\mathcal{L}_XH(x^*)=0.$$

Also,
$$\mathcal{L}_X^2H(x,y)=\nabla \mathcal{L}_XH(x,y)\cdot X(x,y)=f(x,y)g_x(x,y) +g(x,y)g_y(x,y).$$

Substituting $f(x,y)$ and $g(x,y)$ given in \eqref{HamilPert} gives
$$\mathcal{L}_X^2H(x,y)=-(3x^2 + 4 \xi x^3) (-1 + y)^3 +
4 \xi (x^3 + \xi (-1 + x^4 + (-1 + y)^4)) (-1 + y)^3.$$

A direct computation shows that
$$\mathcal{L}_X^2H(x^*)=0.$$

Furthermore,
$$
\begin{array}{rl}
	\mathcal{L}_X^3H(x,y)=&\nabla \mathcal{L}_X^2H(x,y)\cdot X(x,y)\\
	&\\
	=&(x^3 + \xi (-1 + x^4 + (-1 + y)^4)) (-3 (3 x^2 + 4 \xi x^3) (-1 + y)^2 \\
	&\\
	&+12 \xi (x^3 + \xi (-1 + x^4 + (-1 + y)^4)) (-1 + y)^2\\
	&\\
	&+16 \xi^2 (-1 + y)^6) - (-(6 x + 12 \xi x^2) (-1 + y)^3\\
	&\\
	&+4 \xi (3 x^2 + 4 \xi x^3) (-1 + y)^3) (-1 + y)^3,
\end{array}
$$
that provides
$$\mathcal{L}_X^3H(x^*)=0.$$

Finally, computing the fourth Lie derivative at the point $x^*\in\Sigma$ we get
$$\mathcal{L}_X^4H(x^*)=6.$$

So $x^*=(0,0)$ is a regular grazing point of order $4,$ which belongs to the periodic orbit $\Gamma: x^4+(y-1)^4=1 $.

Now we observe that
$$
\begin{array}{ll}
	\mathcal{L}_X^4H(x^*)&=\nabla \mathcal{L}_X^3H(x^*)\cdot X(x^*)\\
	&\\
	&=\left(\dfrac{\partial}{\partial x}\mathcal{L}_X^3H(x^*),\dfrac{\partial}{\partial y}\mathcal{L}_X^3H( x^*)\right)\cdot(f(x^*),g(x^*))\\
	&\\
	&=\dfrac{\partial}{\partial x}\mathcal{L}_X^3H(x^*)f(x^*)+\dfrac{\partial}{\partial y}\mathcal{L}_X ^3H(x^*)g(x^*)\\
	&\\
	&=\dfrac{\partial}{\partial x}\mathcal{L}_X^3H(x^*)f(x^*)+\dfrac{\partial}{\partial y}\mathcal{L}_X ^3H(x^*)\mathcal{L}_XH(x^*)\\
	&\\
	&=\dfrac{\partial}{\partial x}\mathcal{L}_X^3H(x^*)\\
	&\\
	&=6.
\end{array}
$$
This implies that
$$\dfrac{\partial}{\partial x}\mathcal{L}_X^3H(x^*)\neq0.$$ So, from the \textit{Implicit Function Theorem} there is a unique smooth function $x(\epsilon)$ defined in a neighbourhood of $x^*$ such that

$$\mathcal{L}_X^3H(x(\epsilon), -\epsilon)= 0,$$
with $x(0)=0$. Furthermore,
$$x(\epsilon)=x(0)+x'(0)\epsilon+\mathcal{O}(\epsilon^2).$$
By differentiating $\mathcal{L}_X^3H(x(\epsilon), -\epsilon)= 0$ with respect to the $\epsilon$, we have
$$x'(\epsilon)=\dfrac{\dfrac{\partial}{\partial y}\mathcal{L}_X^3H(x(\epsilon), -\epsilon)}{\dfrac{\partial} {\partial x}\mathcal{L}_X^3H(x(\epsilon), -\epsilon)},$$
with
$$
\begin{array}{rl}
	\dfrac{\partial}{\partial y}\mathcal{L}_X^3H(x,y)=&
	(x^3 + \xi (-1 + x^4 + (-1 + y)^4)) (-6 (3 x^2 + 4 \xi x^3) (-1 + y)\\
	& + 24 \xi (x^3 + \xi (-1+ x^4 + (-1 + y)^4)) (-1 + y)+ 144 \xi^2 (-1 + y)^5) \\
	& -3 (-(6 x+ 12 \xi x^2) (-1 + y)^3 + 4 \xi (3 x^2 + 4 \xi x^3) (-1 + y)^3) (-1+ y)^2 \\
	& - (-3 (6 x + 12 \xi x^2) (-1 + y)^2 + 12 \xi (3 x^2 + 4 \xi x^3) (-1+ y)^2) (-1 + y)^3\\
	&+4 \xi (-3 (3 x^2 + 4 \xi x^3) (-1 + y)^2 + 12 \xi (x^ 3+ \xi (-1 + x^4+ (-1 + y)^4)) (-1 + y)^2\\
	&+ 16 \xi^2 (-1 + y)^6) (-1 + y)^3,
\end{array}$$
and
$$
\begin{array}{rl}
	\dfrac{\partial}{\partial x}\mathcal{L}_X^3H(x,y)=&
	(x^3 + \xi (-1 + x^4 + (-1 + y)^4)) (-3 (6 x + 12 \xi x^2) (-1 + y)^2\\
	&+ 12 \xi (3 x^2 + 4 \xi x^3) (-1 + y)^2) + (3 x^2 +
	4 \xi x^3) (-3 (3 x^2 \\
	&+ 4 \xi x^3) (-1 + y)^2 +
	12 \xi (x^3 + \xi (-1 + x^4 \\
	&+ (-1 + y)^4)) (-1 + y)^2 +16 \xi^2 (-1 + y)^6) - (-(6 \\
	& + 24 \xi x) (-1 + y)^3 + 4 \xi (6 x + 12 \xi x^2) (-1 + y)^3) (-1 + y)^3.
\end{array}
$$
So, we obtain
$$x'(0)=\dfrac{\dfrac{\partial}{\partial y}\mathcal{L}_X^3H(x^*)}{\dfrac{\partial}{\partial x}\mathcal {L}_X^3H(x^*)}=\dfrac{-32\xi^3}{3},$$ that implies that

$$x(\epsilon)=\dfrac{-32\xi^3}{3}\epsilon+\mathcal{O}(\epsilon^2).$$

Therefore, given $\mathrm{x_1}=(x(\epsilon),-\epsilon) \in \Pi_3=\{(x,y)\in \mathbb{R}^2;\,\,\mathcal{L }_X^3H(x,y)=0\}$ we have
$$\mathrm{x_1}=\left(\dfrac{-32\xi^3}{3}\epsilon+\mathcal{O}(\epsilon^2),-\epsilon\right).$$

Now, from Theorem \ref{TeoZDM} we get

$$
ZDM(\mathrm{x}_1)=\mathrm{x}_1-\dfrac{(4!)^{3/4}}{3!}W(\mathrm{x}_1)(\mathcal{L} _X^4H(x^*))^{1/4}\epsilon^{3/4}+ \mathcal{O}(\epsilon).$$
So, for $\mathrm{x_1}=\bigg(\dfrac{-32\xi^3}{3}\epsilon+\mathcal{O}(\epsilon^2),-\epsilon\bigg)\in\Pi_3$ it follows that
$$
\begin{array}{ll}
	ZDM(\mathrm{x}_1)=&\bigg(\dfrac{-32\xi^3}{3}\epsilon+\mathcal{O}(\epsilon^2),-\epsilon\bigg)\\
	&\\
	&-\dfrac{(4!)^{3/4}(6)^{1/4}}{3!}W\bigg(\dfrac{-32\xi^3}{3}\epsilon+\mathcal{O} (\epsilon^2),-\epsilon\bigg)\epsilon^{3/4}+ \mathcal{O}(\epsilon).
\end{array}
$$

Therefore,
\begin{equation}\label{ZDMexemplo}
	ZDM(\mathrm{x}_1)=\bigg(\dfrac{-32\xi^3}{3}\epsilon+\mathcal{O}(\epsilon^2),-\epsilon\bigg)-\bigg(\dfrac{4!}{6}\bigg)^{3/4}W\bigg(\dfrac{-32\xi^3}{3}\epsilon+O(\epsilon^2),-\epsilon\bigg )\epsilon^{3/4}+ \mathcal{O}(\epsilon).
\end{equation}

Furthermore, from Theorem \ref{TeoPDM}, we know that
\begin{equation}\label{PDMexamplo}
	PDM(\mathrm{x_1})=\mathrm{x_1}-\Bigg[W(\mathrm{x_1})- \dfrac{\mathcal{L}_W\mathcal{L}_X^3H(\mathrm{ x_1})}{\mathcal{L}_X^4H(x^*)}X(\mathrm{x_1})\Bigg]\dfrac{(4!)^{3/4}}{3!}(\mathcal{L}_X^4H(\mathrm{x ^*}))^{1/4}\epsilon^{3/4}+\mathcal{O}(\epsilon).
\end{equation}

For the particular case where $W(x,y)=\left(k+k_1y, k_2y\right)$,  with $k=cte$ the $ZDM$ becomes
with $k=cste$, we have by \eqref{ZDMexemplo} that
$$
\begin{array}{rl}
ZDM(\mathrm{x}_1)=&\bigg(\dfrac{-32\xi^3}{3}\epsilon+\mathcal{O}(\epsilon^2),-\epsilon\bigg)-\bigg(\dfrac{4!}{6}\bigg)^{3/4}\left(k-k_1\epsilon,-k_2\epsilon\right)\epsilon^{3/4}+ \mathcal{O}(\epsilon).\\
&\\
=&\bigg(\dfrac{-32\xi^3}{3}\epsilon+\mathcal{O}(\epsilon^2),-\epsilon\bigg)-\bigg(\dfrac{4!}{6}\bigg)^{3/4}\left(k,0\right)\epsilon^{3/4}+ \mathcal{O}(\epsilon).
\end{array}
$$

On the other hand, as
$$X(\mathrm{x_1})=X\left(\dfrac{-32\xi^3}{3}\epsilon+\mathcal{O}(\epsilon^2),-\epsilon\right),$$ by using the expression of $X$ given in \eqref{HamilPert} we obtain
$$X(\mathrm{x_1})=\left(-(-1 - \epsilon)^3, -\bigg(\dfrac{32}{3}\bigg)^3 \epsilon^3 \xi^9 -
\xi \left(- (-1 - \epsilon)^4 + \dfrac{(32^4 \epsilon^4 \xi^{12})}{3^4}+\mathcal{O}(\epsilon^6) \right)\right).$$

Using this expression in \eqref{PDMexamplo} a direct computation provides:
$$	PDM(\mathrm{x_1})=\bigg(\dfrac{-32\xi^3}{3}\epsilon+O(\epsilon^2),-\epsilon\bigg)-\dfrac{(4!)^{3/4}(6)^{1/4}}{3!}\left(-k,-\xi\right)\epsilon^{3/4}+ \mathcal{O}(\epsilon).$$

\end{section}

\section*{Appendix}\label{s2}
In this section we briefly describe the basic result of perturbation theory that we need
to obtain Proposition \ref{propZDM}.

Next result guarantees sufficient conditions to estimate the root of a perturbed function. For a complete proof see chapter 3 of \cite{reference9}

\begin{teo}(\cite{reference9})\label{ObsRaiz}
	Consider $f \in C^{1}(\mathbb{R})$ and $g\in C^{m}(\mathbb{R})$. Suppose that $g$ has a root $\overline{\delta}$ of multiplicity m. If g is perturbed by $\epsilon f,$ with $\epsilon\sim 0^+$, then the root of $\epsilon f(x)+ g(x)$ is of the form
	$$\delta=\overline{\delta}+\left(-\epsilon\dfrac{ m!f(\overline{\delta})}{g^{(m)}(\overline{\delta})}\right )^{\frac{1}{m}}e^{\frac{2k\pi i}{m}},
	$$
	with $k\le0<m$. In particular, when $m=1$, we have $k=0$ and
	$$
	\begin{array}{rl}
		\delta=&\overline{\delta}+\left(-\dfrac{\epsilon m!f(\overline{\delta})}{g_{x}(\overline{\delta})}\right) ^{\frac{1}{m}}e^{\frac{2k\pi i}{m}}\\
		&\\
		=&\overline{\delta}-\epsilon\dfrac{f(\overline{\delta})}{g_{x}(\overline{\delta})}.
	\end{array}	
	$$
\end{teo}

\bibliographystyle{plain}

g=\end{document}